\documentclass[aos,seceqn,dvips]{arximspdf}
\usepackage{graphicx}

\doi{10.1214/10-AOS841}
\volume{39}
\issue{1}
\pubyear{2011}
\firstpage{514}
\lastpage{555}

\makeatletter
\newtheorem{theorem}{Theorem}
\newtheorem{prop}{Proposition}
\newtheorem{lemma}{Lemma}
\newtheorem{cor}{Corollary}
\newproclaim{remark}{Remark}
\newproclaim{definition}{Definition}
\makeatother

\begin{document}
\begin{frontmatter}

\title{Wishart distributions for decomposable covariance graph models}
\runtitle{Wishart distributions}

\begin{aug}
\author[A]{\fnms{Kshitij} \snm{Khare}\thanksref{t1}\ead[label=e1]{kdkhare@stat.ufl.edu}}
\and
\author[B]{\fnms{Bala} \snm{Rajaratnam}\thanksref{t2}\corref{}\ead[label=e2]{brajarat@stanford.edu}}
\thankstext{t1}{Supported in part by the B. C. and E. J. Eaves Stanford graduate fellowship.}
\thankstext{t2}{Supported in part by NSF Grants DMS-05-05303, DMS-09-06392 and Grant SUFSC08-SUSHSTF09-SMSCVISG0906.}
\runauthor{K. Khare and B. Rajaratnam}
\affiliation{University of Florida and Stanford University}
\address[A]{Department of Statistics\\
University of Florida\\
102 Griffin-Floyd Hall\\
Gainesville, Florida 32606\\
USA\\
\printead{e1}} 
\address[B]{Department of Statistics\\
 Stanford University\\
Sequoia Hall\\
390 Serra Mall\\
Stanford, California 94305-4065\\
USA\\
\printead{e2}}
\end{aug}

\received{\smonth{7} \syear{2009}}
\revised{\smonth{6} \syear{2010}}

\begin{abstract}

Gaussian covariance graph models encode marginal independence among the
components of a multivariate random vector by means of a graph $G$.
These models are distinctly different from the traditional
concentration graph models (often also referred to as Gaussian
graphical models or covariance selection models) since  the zeros in
the parameter are now reflected in the covariance matrix $\Sigma$, as
compared to the concentration matrix $\Omega = \Sigma^{-1}$. The
parameter space of interest for covariance graph models is the cone
$P_G$ of positive definite matrices with fixed zeros corresponding to
the missing edges of $G$. As in Letac and Massam [\textit{Ann.
Statist.} \textbf{35} (2007) 1278--1323], we consider the case where
$G$ is decomposable. In this paper, we construct on the cone $P_G$ a
family of Wishart distributions which\ serve a similar purpose in the
covariance graph setting as those constructed by Letac and Massam
[\textit{Ann. Statist.} \textbf{35} (2007) 1278--1323] and Dawid and
Lauritzen [\textit{Ann. Statist.} \textbf{21} (1993)  1272--1317] do in
the concentration graph setting. We proceed to undertake a rigorous
study of these ``covariance'' Wishart distributions and derive several
deep and useful properties of this class. First, they form a rich
conjugate family of priors with multiple shape parameters for
covariance graph models. Second, we show how to sample from these
distributions by using a block Gibbs sampling algorithm and prove
convergence of this block Gibbs sampler. Development of this class of
distributions enables Bayesian inference, which, in turn, allows for
the estimation of $\Sigma$, even in the case when the sample size is
less than the dimension of the data (i.e., when ``$n < p$''), otherwise
not generally possible in the maximum likelihood framework. Third, we
prove that when $G$ is a homogeneous graph, our covariance priors
correspond to standard conjugate priors for appropriate directed
acyclic graph (DAG) models. This correspondence enables closed form
expressions for normalizing constants and expected values, and also
establishes hyper-Markov properties for our class of priors. We also
note that when $G$ is homogeneous, the family $\mathit{IW}_{Q_G}$ of Letac and
Massam [\textit{Ann. Statist.} \textbf{35} (2007) 1278--1323] is a
special case of our covariance Wishart distributions. Fourth, and
finally, we illustrate the use of our family of conjugate priors on
real and simulated data.
\end{abstract}



\end{frontmatter}

\section{Introduction}\label{sec1}

Due to recent advances in science and information technology, there has
been a huge influx of high-dimensional data from various fields such as
genomics, environmental sciences, finance and the social sciences.
Making sense of all the many complex relationships and multivariate
dependencies present in the data, formulating correct models and
developing inferential procedures is one of the major challenges in
modern day statistics. In parametric models, the covariance or
correlation matrix (or its inverse) is the fundamental object that
quantifies relationships between random variables. Estimating the
covariance matrix in a sparse way is crucial in high-dimensional
problems and enables the detection of the most important relationships.
In this light, graphical models have served as tools to discover
structure in high-dimensional data.

 The primary aim of this paper is to develop a new family of
conjugate prior distributions for covariance graph models (a subclass
of graphical models) and study the properties of this family of
distributions. It is shown in this paper that these properties are
highly attractive for Bayesian inference in high-dimensional settings.
In covariance graph models, specific entries of the covariance matrix
are restricted to be zero, which implies marginal independence in the
Gaussian case. Covariance graph models correspond to curved exponential
families and are distinctly different from the well-studied
concentration graph models, which, in turn, correspond to natural
exponential families.

 A rich framework for Bayesian inference for natural exponential
families has been established in the last three decades, starting with
the seminal and celebrated work of Diaconis and Ylvisaker
\cite{dicnsylskr} that laid the foundations for constructing conjugate
prior distributions for natural exponential family models. The
Diaconis--Ylvisaker (henceforth referred to as ``DY'') conjugate priors
are characterized by posterior linearity of the mean. An analogous
framework for curved exponential families is not available in the
literature.

 Concentration graph models (or covariance selection models)
were one of the first graphical models to be formally introduced to the
statistics community. These models reflect conditional independencies
in multivariate probability distributions by means of a graph. In the
Gaussian case, they induce sparsity or zeros in the inverse covariance
matrix and correspond to natural exponential families. In their
pioneering work, Dawid and Lauritzen \cite{dawidlrtzn} developed the DY
prior for this class of models. In particular, they introduced the
hyper-inverse Wishart as the DY conjugate prior for concentration graph
models. In a recent major contribution to this field, a rich family of
conjugate priors that subsumes the DY class has been developed by Letac
and Massam \cite{ltcmssmwdg}. Both the hyper-inverse Wishart priors and
the ``Letac--Massam'' priors have attractive properties which enable
Bayesian inference, with the latter allowing multiple shape parameters
and hence being suitable in high-dimensional settings. Bayesian
procedures corresponding to these Letac--Massam priors have been
derived in a decision theoretic framework in the recent work of
Rajaratnam, Massam and Carvalho \cite{rajmasscar}.
\setcounter{footnote}{2}

 Consider an undirected\footnote{We shall use dotted edges for
our graphs, in keeping  with the notation in the literature;
bi-directed edges have also been used for representing covariance
graphs.} graph $G$ with a finite set of  vertices $V$ (of size $p$) and
a finite set $E$ of edges between these vertices, that is, $G = (V,E)$.
The Gaussian covariance graph model corresponding to the graph $G$ is
the collection of $p$-variate Gaussian distributions with covariance
matrix $\Sigma$ such that $\Sigma_{ij} = 0$ whenever $(i,j) \notin E$.
This class of models was first formally introduced by Cox and Wermuth
\cite{cxwrmthcvg,cxwrmthmai}. In the frequentist setting, maximum
likelihood estimation in covariance graph models has been a topic of
interest in recent years. Many iterative methods that obtain the
maximum likelihood estimate have been proposed in the literature. The
graphical modeling software MIM in Edwards \cite{edwardsigm} fits these
models by using the ``dual likelihood method'' from Kauermann
\cite{krmnndlldm}. In Wermuth, Cox and Marchetti \cite{wrmthcxmrt}, the authors
derive asymptotically efficient approximations to the maximum
likelihood estimate in covariance graph models for exponential
families. Chaudhuri, Drton and Richardson \cite{chaudrtric} propose an iterative
conditional fitting algorithm for maximum likelihood estimation in this
class of models. Covariance graph models have also been used in
applications in Butte et al. \cite{butaslgoko}, Grzebyk, Wild and
Chouaniere \cite{grzwildcho}, Mao, Kschischang and Frey \cite{maokscfrey} and others.

 Although Gaussian covariance graph models are simple and
intuitive to understand, no comprehensive theoretical framework for
Bayesian inference for this class of models has been developed in the
literature. In that sense, Bayesian inference for covariance graph
models has been an open problem since the introduction of these models
by Cox and Wermuth \cite{cxwrmthcvg,cxwrmthmai} more than fifteen years
ago. The main difficulty is that these models give rise to curved
exponential families. The zero restrictions on the entries of the
covariance matrix $\Sigma$ translate into complicated restrictions on
the corresponding entries of the natural parameter, $\Omega =
\Sigma^{-1}$. Hence, the sparseness in $\Sigma$ does not translate into
sparseness in $\Sigma^{-1}$ and thus a covariance graph model cannot be
viewed as a concentration graph model. No general theory is available
for Bayesian inference in curved exponential families for continuous
random variables, akin to the Diaconis--Ylvisaker \cite{dicnsylskr} or
standard conjugate theory for natural exponential families.

 There are several desirable properties that one might want when
constructing a class of priors, but one of the foremost requirements is
to be able to compute quantities such as the mean or mode of the
posterior distribution, either in closed form or by sampling from the
posterior distribution by a simple mechanism. This is especially
important in high-dimensional situations, where computations are
complex and can become infeasible very quickly. Another desirable and
related feature is conjugacy, that is, the class of priors is such that
the posterior distribution also belongs to this class. Among other
things, this increases the prospects of obtaining closed form Bayes
estimators and can also add to the interpretability of the
hyper-parameters. The class of Wishart distributions developed by Letac
and Massam \cite{ltcmssmwdg} (and later used for flexible Bayesian
inference for concentration graph models by Rajaratnam, Massam and Carvalho
\cite{rajmasscar}), known as the $\mathit{IW}_{P_G}$ family of distributions,
are not appropriate for the covariance graph setting. There is the
additional option of using  the $\mathit{IW}_{Q_G}$ class as priors for this
situation. We, however, establish that the posterior distribution fails
to belong to the same class and there are no known results for
computing the posterior mean or mode, either in closed form or by
sampling from the posterior distribution.

 A principal objective of this paper is to develop a framework
for Bayesian inference for Gaussian covariance graph models. We proceed
to construct a rich and flexible class of conjugate Wishart
distributions, with multiple shape parameters, on the space of positive
definite matrices with fixed zeros, that corresponds to a decomposable
graph $G$. This class of distributions is specified up to a normalizing
constant, and conditions under which this normalizing constant can be
evaluated in closed form are derived. We explore the distributional
properties of our class of priors and, in particular, show that the
parameter can be partitioned into blocks so that the conditional
distribution of each block, given the others, is tractable. Based on
this property, we propose a block Gibbs sampling algorithm to simulate
from the posterior distribution. We proceed to formally prove the
convergence of this block Gibbs sampler. Our priors yield proper
inferential procedures, even in the case when the sample size $n$ is
less than the dimension $p$ of the data, whereas maximum likelihood
estimation is, in general, only possible when $n \geq p$ (in fact, in
the homogeneous case, it can be shown that the condition $n \geq p$ is
actually also necessary, thus highlighting the fact that results from
the concentration graph setting do not carry over to the covariance
model setting). We also show that our covariance Wishart distributions
are, in the decomposable nonhomogeneous case, very different from the
Letac--Massam priors $W_{P_G}$ and $\mathit{IW}_{Q_G}$. However, when the
underlying graph $G$ is homogeneous, the Letac--Massam $\mathit{IW}_{Q_G}$
 priors are a special case of our distributions. We establish, in the
 homogeneous setting, a correspondence between the covariance priors
 in this paper and the natural conjugate priors for appropriate
 directed acyclic graph (DAG) models. This correspondence helps us
 to explicitly evaluate quantities like the normalizing constant
 and the posterior mean of the covariance matrix in closed form.
 In this scenario, we also show that our class of priors satisfies
 the strong directed hyper-Markov property (as introduced in Dawid and
 Lauritzen \cite{dawidlrtzn} for concentration graph models).
 It should be pointed out that these aforementioned results for
 homogeneous graphs can also be established directly, without
 exploiting the correspondence with the DAG models. The direct
 approach is self-contained, whereas the latter invokes an external
 result which states that for the restrictive class of homogeneous
 graphs, covariance graph models and DAGs are Markov equivalent.

We noted above that for concentration graph models or the traditional
Gaussian graphical models, a rich theory has been established by Dawid
and Lauritzen \cite{dawidlrtzn}, who derive the single parameter DY
conjugate prior for these models, and by Letac and Massam
\cite{ltcmssmwdg}, who derive a larger flexible class with multiple
shape parameters. In essence, this paper is the analog of the results
in the two aforementioned papers in the covariance graph model setting,
with parallel results, all of which are contained in a single
comprehensive piece. Hence, this work completes the powerful theory
that has been developed in the mathematical statistics literature for
decomposable models.

 We also point out that a class of priors in the recent work
\cite{ghrmnslvgs} is a special case of our class of flexible covariance
Wishart distributions.\footnote{This is in a similar spirit to the way
in which the HIW prior of Dawid and Lauritzen \cite{dawidlrtzn} is a
special case of the generalized family of Wishart distributions
proposed by Letac and Massam \cite{ltcmssmwdg} for the concentration
graph setting.} Our family allows multiple shape parameters, as
compared to a single shape parameter, and hence yields a richer class
suitable to high-dimensional problems. Moreover, we show that their
iterative algorithm to sample from the posterior is different from
ours. Since the authors do not undertake a theoretical investigation of
the convergence properties of their algorithm, it is not clear if it
does indeed converge to the desired distribution. On the other hand, we
proceed to formally prove that our algorithm converges to the desired
distribution. The remaining sections of this paper are considerably
different from \cite{ghrmnslvgs} since we undertake a rigorous
probabilistic analysis of our conjugate Wishart distributions for
covariance graph models, whereas they give a useful and novel treatment
of latent variables and mixed graph models in a machine learning
context.

 This paper is structured as follows. Section \ref{prelimmtrl}
introduces the required preliminaries and notation. In Section
\ref{priorcnsrt}, the class of covariance Wishart distributions is
formally constructed. Conjugacy to the class of covariance graph models
and sufficient conditions for integrability are established. Comparison
with the Letac--Massam $\mathit{IW}_{Q_G}$ priors, which are not, in general,
conjugate in the covariance graph setting, is also undertaken. In
Section \ref{cndtlcmptn}, a block Gibbs sampler which enables sampling
from the posterior distribution is proposed and the corresponding
conditional distributions are derived. Thereafter, a formal proof of
convergence of this block Gibbs sampler is provided. In Section
\ref{hprmrkvpty}, we restrict ourselves to the case when $G$ is a
homogeneous graph. We examine the distributional properties of our
class of priors in this section and prove that the covariance priors
introduced in this paper correspond to natural conjugate priors for DAG
models in the homogeneous setting. This correspondence helps in
establishing closed form expressions for normalizing constants,
expected values and hyper-Markov properties for our class of priors for
$G$ homogeneous. Finally, we illustrate the use of our family of
conjugate priors and the methodology developed in this paper on a real
example, as well as on simulated data. The \hyperref[app]{Appendix} contains the proofs
of some of the results stated in the main text.

\section{Preliminaries} \label{prelimmtrl}
In this section, we give the necessary notation, background and
preliminaries that are needed in subsequent sections.

\subsection{Modified Cholesky decomposition}
If $\Sigma$ is a positive definite matrix, then there exists a unique
decomposition
\begin{equation} \label{inverselwr}
\Sigma = LDL^T,
\end{equation}
where $L$ is a lower-triangular matrix with diagonal entries equal to 1
and $D$ a diagonal matrix with positive diagonal entries. This
decomposition of $\Sigma$ is referred to as the \textit{modified Cholesky
decomposition} of $\Sigma$ (see \cite{pourahmadi}). We now provide a
formula that explicitly computes the inverse of a lower-triangular
matrix with $1$'s on the diagonal, such as those that appear in
(\ref{inverselwr}).

\begin{prop} \label{inverselrt}
Let $L$ be an $m \times m$ lower-triangular matrix with $1$'s on the
diagonal. Let
\[
\mathcal{A} = \bigcup_{r=2}^m \bigl\{\bolds{\tau}\dvtx \bolds{\tau} \in \{1, 2,
\ldots, m\}^r, \tau_i < \tau_{i-1}\ \forall  2 \leq i \leq r \bigr\}
\]
and
\[
L_{\bolds{\tau}} = \prod_{i=2}^{\operatorname{dim}(\bolds{\tau})} L_{\tau_{i-1} \tau_i}
 \qquad \forall \bolds{\tau} \in \mathcal{A},
\]
where $\operatorname{dim}(\tau)$ denotes the length of the vector $\tau$. Then,
$L^{-1} = N$, where $N$ is lower-triangular matrix with $1$'s on the
diagonal and, for $i>j$,
\[
N_{ij}=\sum_{\tau\in {\mathcal A},\tau_1=i,\tau_{\operatorname{dim}{\tau}}=j}(-1)^{\operatorname{dim}(\tau)-1}
\prod_{i=2}^{\operatorname{dim}(\tau)}L_{\tau_{i-1}\tau_i}.
\]
\end{prop}

The proof is provided in the \hyperref[app]{Appendix}.

 An undirected graph $G$ is a pair $(V,E)$, where $V$ is a
permutation\footnote{The ordering in $V$ is emphasized here since the
elements of $V$ will later correspond to rows or columns of matrices.}
of the set $\{1,2, \ldots, m\}$ denoting the set of vertices of $G$.
The set $E \subseteq V \times V$ denotes the set of edges in the graph.
If vertices $u$ and $v$ are such that $(u,v) \in E,$ then we say that
there is an edge between $u$ and $v$. It is also understood that $(u,v)
\in E$ implies that $(v,u) \in E$, that is, the edges are undirected.
Although the dependence of $G = (V,E)$ on the particular ordering in
$V$ is often suppressed, the reader should bear in mind that unlike
traditional graphs, the graphs defined above are not equivalent up to
permutation of the vertices\footnote{This has been done for notational
convenience, as will be seen later.} modulo the edge structure. Below,
we describe  two classes of graphs which play a central role in this
paper.

\subsection{Decomposable graphs} \label{dmpgrphpty}

An undirected graph $G$ is said to be \textit{decomposable} if any induced
subgraph does not contain a cycle of length greater than or equal to
four. The reader is referred to Lauritzen \cite{lrtzngphmd} for all of
the common notions of graphical models (and, in particular,
decomposable graphs) that we will use here. One such important notion
is that of a perfect order of the cliques. Every decomposable graph
admits a  perfect order of its cliques. Let $(C_1, C_2, \ldots, C_k)$
be one such perfect order of the cliques of the graph $G$. The {\it
history} for the graph is given by $H_1 = C_1$ and
\[
H_j = C_1 \cup C_2 \cup \cdots \cup C_j,\qquad j = 2, 3, \ldots, k,
\]
and the \textit{minimal separators} of the graph are given by
\[
S_j = H_{j-1} \cap C_j,\qquad j = 2, 3, \ldots, k.
\]
Let
\[
R_j = C_j\setminus H_{j-1} \qquad\mbox{for } j = 2, 3, \ldots, k.
\]

Let $k' \leq k - 1$ denote the number of distinct separators and $\nu
(S)$ denote the multiplicity of $S$, that is, the number of $j$ such
that $S_j = S$. Generally, we will denote by $\mathcal{C}$ the set of
cliques of a graph and by $\mathcal{S}$ its set of separators.

Now, let $\Sigma$ be an arbitrary positive definite matrix with
zero restrictions according to $G = (V,E)$,\footnote{It is emphasized
here that the ordering of the vertices reflected in $V$ plays a crucial
role in the definitions and results that follow.} that is,
$\Sigma_{ij} = 0$ whenever $(i,j) \notin E$. It is known that if $G$ is
decomposable, then there exists an ordering of the vertices such that
if $\Sigma = LDL^T$ is the modified Cholesky decomposition
corresponding to this ordering, then, for $i > j$,
\begin{equation} \label{choleskirn}
L_{ij} = 0 \qquad\mbox{whenever } (i,j) \notin E.
\end{equation}

Although the ordering is not unique in general, the existence of such
an ordering characterizes decomposable graphs (see \cite{plsnpwrsmh}).
A constructive way to obtain such an ordering is given as follows.
Label the vertices in descending order, starting with vertices in $C_1,
R_2, R_3, \ldots, R_k$, with vertices belonging to a particular set
being ordered arbitrarily (see \cite{lrtzngphmd,plsnpwrsmh,wermuthcsp}
for more details).

\subsection{The spaces $P_G$, $Q_G$ and $\mathcal{L}_G$} \label{spaceintro}

An $m$-dimensional Gaussian covariance graph model\footnote{A brief
overview of the literature in this area is provided in the
\hyperref[sec1]{Introduction}.} can be represented by the class of multivariate normal
distributions with fixed zeros in the covariance parameter (i.e.,
marginal independencies) described by a given graph $G = (V,E)$. That
is, if $(i,j)\notin E$, then the $i$th and $j$th components of the
multivariate random vector are marginally independent. Without loss of
generality, we can assume that these models have mean zero and are
characterized by the parameter set $P_G$ of positive definite
covariance matrices $\Sigma$ such that $\Sigma_{ij}=0$ whenever the
edge $(i,j)$ is not in $E$. Following the notation in
\cite{ltcmssmwdg,rajmasscar} for $G$ decomposable, we define $Q_G$ to
be the space on which the free elements of the precision matrices (or
inverse covariance matrices) $\Omega$ live.

More formally, let $M$ denote the set of symmetric matrices of order
$m$,  $M^+_m\subset M$ the cone of positive definite matrices
(abbreviated as ``$>0$''), $I_G$ the linear space of symmetric
incomplete matrices $x$ with missing entries $x_{ij}, (i,j) \notin E$,
and $\kappa\dvtx M \mapsto I_G$ the projection of $M$ into $I_G$. The
parameter set of the precision matrices of Gaussian covariance graph
models can also be described as the set of incomplete matrices $\Omega=
\kappa(\Sigma^{-1}), \Sigma\in P_G$. The entries $\Omega_{ij}, (i,j)
\notin E,$ are not free parameters of the precision matrix for Gaussian
covariance graph models (see \cite{ltcmssmwdg,rajmasscar} for details).
We are therefore led to consider the two cones
\begin{eqnarray}\label{pg}
P_G&=&\{y\in M^+_m| y_{ij}=0, (i,j)\notin E\},\\\label{qg}
Q_G&=&\{x\in I_G| x_{C_i}>0, i=1,\ldots,k\},
\end{eqnarray}
where $P_G\subset Z_G$, $Q_G\subset I_G$ and  $Z_G$ denotes the linear
space of symmetric matrices with zero entries $y_{ij}, (i,j)\notin
E$. Furthermore Grone et al. \cite{Grone84} prove that for $G$
decomposable, the spaces $P_G$ and $Q_G$ are isomorphic (once more, see
\cite{ltcmssmwdg,rajmasscar} for details).

We now introduce new spaces $\mathcal{L}_G$ and ${\bolds
\Theta}_G$ (the modified Cholesky space) that will be needed in our
subsequent analysis\footnote{These spaces are not defined in
\cite{ltcmssmwdg,rajmasscar}.}:
\begin{eqnarray*}
\mathcal{L}_G &=& \{L\dvtx   L_{ij} = 0 \mbox{ whenever } i < j, \mbox{ or } (i,j) \notin E,
\mbox{ and } L_{ii} = 1, \forall  1 \leq i,j \leq m\};\\
{\bolds\Theta}_G &=& \{\bolds{\theta} = (L,D)\dvtx   L \in \mathcal{L}_G,
 D \mbox{ diagonal with } D_{ii} > 0\ \forall 1 \leq i \leq m\}.
\end{eqnarray*}
We define the mapping $\psi\dvtx {\bolds\Theta}_G \rightarrow M_m^+$ as
follows:
\begin{equation} \label{bijecthpcg}
\psi (L,D) = LDL^T.
\end{equation}
This mapping $\psi$ plays an important role in our analysis and shall
be studied later.

\subsection{Homogeneous graphs} \label{hmgnsgphpy}
A graph $G = (V,E)$ is defined to be homogeneous if, for all $(i,j)\in
E$, either
\[
\{u\dvtx  u = j \mbox{ or } (u,j) \in E\} \subseteq \{u\dvtx  u = i \mbox{ or }
(u,i) \in E\}
\]
or
\[
\{u\dvtx  u = i \mbox{ or } (u,i) \in E\} \subseteq \{u\dvtx  u = j \mbox{ or } (u,j) \in E\}.
\]
Equivalently, a graph $G$ is said to be homogeneous if it is
decomposable and does not contain the graph
$\stackrel{1}{\bullet} - \stackrel{2}{\bullet} -
\stackrel{3}{\bullet} - \stackrel{4}{\bullet}$, denoted
by $A_4$, as an induced subgraph. Homogeneous graphs have an equivalent
representation in terms of directed rooted trees, called \textit{Hasse
diagrams}. The reader is referred to \cite{ltcmssmwdg} for a detailed
account of the properties of homogeneous graphs. We write $i
\rightarrow j$ whenever
\[
\{u\dvtx  u = j \mbox{ or } (u,j) \in E\} \subseteq \{u\dvtx  u = i \mbox{ or }
(u,i) \in E\}.
\]
Denote by $R$ the equivalence relation on $V$ defined by
\[
i R j \quad\Leftrightarrow\quad i \rightarrow j \mbox{ and } j \rightarrow i.
\]

 Let $\bar{i}$ denote the equivalence class in $V/R$ containing
$i$. The Hasse diagram of $G$ is defined as a directed graph with
vertex set $V_H = V/R = \{\bar{i}\dvtx i \in V\}$ and edge set $E_H$
consisting of directed edges with $(\bar{i}, \bar{j}) \in E_H$ for
$\bar{i} \neq \bar{j}$ if the following holds: $i \rightarrow j$ and
$\nexists k$ such that $i \rightarrow k \rightarrow j$, $\bar{k} \neq
\bar{i}$, $\bar{k} \neq \bar{j}$.

\begin{figure}[b]
\begin{tabular}{@{}cc@{}}

\includegraphics{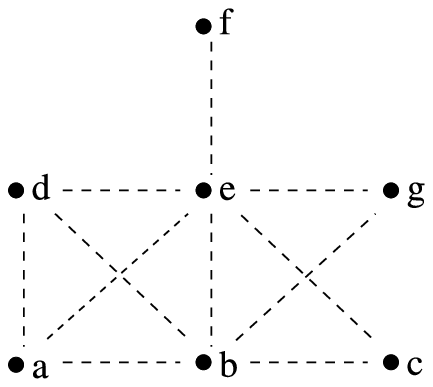}
&\includegraphics{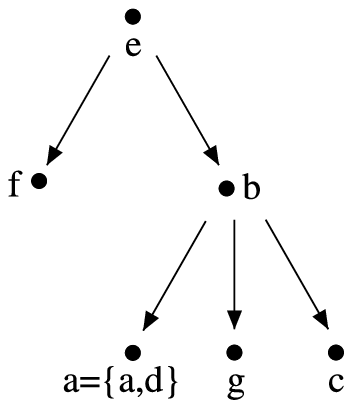}\\
(a)&(b)
\end{tabular}
\caption{\textup{(a)} An example of a homogeneous graph with $7$ vertices;
\textup{(b)}
the corresponding Hasse diagram.} \label{hassegraph}
\end{figure}

 If $G$ is a homogeneous graph, then the Hasse diagram described
above is a directed rooted tree such that the number of children of a
vertex is never equal to one. It was proven in \cite{ltcmssmwdg} that
there is a one-to-one correspondence between the set of homogeneous
graphs and the set of directed rooted trees with vertices weighted by
positive integers [$w(\bar{i}) = |\bar{i}|$], such that no vertex has
exactly one child. Also, when $i R j$, we say that $i$ and $j$ are
\textit{twins} in the Hasse diagram of $G$. Figure \ref{hassegraph}
provides an example of a homogeneous graph with seven vertices and the
corresponding Hasse diagram.

The following proposition for homogeneous graphs plays an
important role in our analysis.

\begin{prop} \label{hmgnsrstrt}
If $G$ is a homogeneous graph, then there exists an ordering of the
vertices, such that, for this ordering:
\begin{enumerate}
\item $\Sigma \in P_G \Leftrightarrow L \in \mathcal{L}_G$, where
$\Sigma = LDL^T$ is the modified Cholesky decomposition of $\Sigma$;
\item $L \in \mathcal{L}_G \Leftrightarrow L^{-1} \in \mathcal{L}_G$.
\end{enumerate}
\end{prop}

The proof of this proposition is well known and so is omitted for the
sake of brevity (see \cite{adsnwgnrwd,khrrjtmwdg,roverchldn}). We now
describe a procedure for ordering the vertices, under which Proposition
\ref{hmgnsrstrt} holds. Given a homogeneous graph $G$, we first
construct the Hasse diagram for $G$. The vertices are labeled in
descending order, starting from the root of the tree. If the
equivalence class at any node has more than one element, then they are
labeled in any order. Hereafter, we shall refer to this ordering scheme
as the \textit{Hasse perfect vertex elimination scheme}. For example, if
we apply this ordering procedure to the graph in Figure
\ref{hassegraph}, then the resulting labels are $\{a,b,c,d,e,f,g\}
\rightarrow \{4,5,1,3,7,6,2\}$.

\subsection{Vertex ordering}

Let $G = (V,E)$ be an undirected decomposable graph with vertex set $V
= \{1,2, \ldots, m\}$ and edge set $E$. Let $S_{V}$ denote the
permutation group associated with $V$. For any $\sigma \in S_{V}$, let
$G_\sigma := (\sigma (V), E_\sigma)$, where $(u,v) \in E_\sigma$ if and
only if $(\sigma^{-1} (u), \sigma^{-1} (v)) \in E$. Let $S_D \subseteq
S_{V}$ denote the subset of permutations $\sigma$ of $V$ such that, for
any $\Sigma \in M_m^+$ with $\Sigma = LDL^T$, $L \in
\mathcal{L}_{G_\sigma} \Leftrightarrow \Sigma \in P_{G_\sigma}$. Hence,
for every $\sigma \in S_D$, the  mapping $\psi_\sigma \dvtx  {\bolds
\Theta}_{G_\sigma} \rightarrow M_m^+$ defined in (\ref{bijecthpcg}) is
a bijection from ${\bolds\Theta}_{G_\sigma}$ to $P_{G_\sigma}$. In
particular, the ordering corresponding to any perfect vertex
elimination scheme lies in $S_D$ (see Section \ref{dmpgrphpty}). If $G$
is homogeneous, let $S_H \subseteq S_D$ denote the subset of
permutations $\sigma$ of $V$ such that $L \in \mathcal{L}_{G_\sigma}
\Leftrightarrow L^{-1} \in \mathcal{L}_{G_\sigma}$. In particular, any
ordering of the vertices corresponding to the Hasse perfect vertex
elimination scheme lies in $S_H$ (see Section \ref{hmgnsgphpy}). The
above defines a nested triplet of permutations of $V$ given by $S_H
\subseteq S_D \subseteq S_{V}$.

\section{Wishart distributions for covariance graphs} \label{priorcnsrt}

Let $G = (V,E)$ be an undirected decomposable graph with vertex set $V$
and edge set $E$. We assume that the vertices in $V$ are ordered so
that $V \in S_D$. The covariance graph model associated with $G$ is the
family of distributions
\begin{eqnarray*}
\mathcal{G}
&=& \{ \mathcal{N}_m ({\mathbf 0}, \Sigma)\dvtx  \Sigma \in P_G \}\\
&\cong& \{ \mathcal{N}_m ({\mathbf 0}, LDL^T)\dvtx  (L,D) \in {\bolds
\Theta}_G \}.
\end{eqnarray*}

Consider the class of measures on ${\bolds \Theta}_G$ with density
[with respect to\break $\prod_{i>j, (i,j) \in E}dL_{ij}\,
\prod_{i=1}^m dD_{ii}$]
\begin{equation} \label{flexiblecjg}
\hspace*{15pt}\widetilde{\pi}_{U, \bolds{\alpha}} (L, D) = e^{- (\operatorname{tr} (
(LDL^T)^{-1} U ) + \sum_{i=1}^m \alpha_i \log D_{ii} )/2},\qquad \theta =
(L, D) \in {\bolds \Theta}_G.
\end{equation}
These measures are parameterized by a positive definite matrix $U$ and
a vector $\bolds{\alpha} \in \mathbb{R}^m$ with nonnegative entries.
Let us first establish some notation:
\begin{itemize}
\item $\mathcal{N} (i) := \{j\dvtx  (i,j) \in E\};$
\item $\mathcal{N}^\prec (i) := \{j\dvtx  (i,j) \in E, i > j\};$
\item $U^{\prec i} := ((U_{kl}))_{k,l \in \mathcal{N}^\prec (i)};$
\item $U^{\preceq i}:= ((U_{kl}))_{k,l \in \mathcal{N}^\prec (i) \cup \{i\}};$
\item $U^\prec_{\cdot i} := (U_{ki})_{k \in \mathcal{N}^\prec (i)}.$
\end{itemize}
Let
\[
z_G (U, \bolds{\alpha}) := \int e^{- ( \operatorname{tr} ( (LDL^T)^{-1} U ) +
\sum_{i=1}^m \alpha_i \log D_{ii} )/2}\,dL\,dD.
\]
If $z_G (U,\bolds{\alpha}) < \infty$, then $\widetilde{\pi}_{U,\bolds{\alpha}}$ can be normalized
to obtain a probability measure. A~sufficient condition for the existence of a normalizing constant for
$\widetilde{\pi}_{U, \bolds{\alpha}} (L,D)$ is provided in the
following proposition.

\begin{theorem} \label{nrmcnstdml}
Let $dL := \prod_{(i,j) \in E, i>j} dL_{ij}$ and $dD := \prod_{i=1}^m
dD_{ii}$. Then,
\[
\int_{{\bolds \Theta}_G} e^{- ( \operatorname{tr} ( (LDL^T)^{-1} U ) +
\sum_{i=1}^m \alpha_i \log D_{ii} )/2}\,dL\,dD < \infty
\]
if
\[
\alpha_i > |\mathcal{N}^\prec (i)| + 2\qquad\forall  i = 1,2, \ldots, m.
\]
\end{theorem}

As the proof of this proposition is rather long and technical, it is
deferred to the \hyperref[app]{Appendix}. The normalizing constant
$z_G (U,\bolds{\alpha})$ is not generally available in closed form. Let us consider a
simple example to illustrate the difficulty of computing the
normalizing constant explicitly.

Let $G = A_4$, that is, the path on four vertices, or
$\stackrel{1}{\bullet} - \stackrel{2}{\bullet} -
\stackrel{3}{\bullet} - \stackrel{4}{\bullet}$. Note that
this is a decomposable (but not homogeneous) graph. The restrictions on
$L$ are $L_{31} = L_{41} = L_{42} = 0$. Let $U \in P_G$ and $
\bolds{\alpha} = (4,4,4,4)$. Then, after integrating out the elements $D_{ii},
1 \leq i \leq 4$ (recognizing them as inverse-gamma integrals) and
transforming the entries of $L$ to the independent entries of $L^{-1}$
(as in the proof of Proposition \ref{nrmcnstdml}), the normalizing
constant reduces to an integral of the form
\begin{eqnarray*}
\int_{\mathbb{R}^3}&&(U_{22} + 2 U_{12}x_1 +
U_{11}x_1^2)^{-1}\\
&&{}\times(U_{11} x_1^2 x_2^2 +
U_{22} x_2^2 + U_{33} + 2 U_{12} x_1 x_2^2  + 2 U_{13} x_1 x_2 + 2 U_{23} x_2)^{-1}\\
&&{}\times(U_{11} x_1^2 x_2^2 x_3^2 + U_{22} x_2^2 x_3^2 +
U_{33} x_3^2 + U_{44}+ 2 U_{12} x_1 x_2^2 x_3^2\\
&&\quad\hphantom{(}{} + 2 U_{13} x_1 x_2
x_3^2 + 2 U_{14} x_1 x_2 x_3+ 2 U_{23} x_2 x_3^2 + 2 U_{24} x_2 x_3 +
U_{34} x_3)^{-1}\,d{\mathbf x}.
\end{eqnarray*}
The above integral does not seem to be computable by standard
techniques for general $U$. Despite this inherent difficulty, we
propose a novel method which allows sampling from this rich family of
distributions (see Section \ref{cndtlcmptn}).

We will show later that the condition in Theorem \ref{nrmcnstdml} is
necessary and sufficient for the existence of a normalizing constant
for homogeneous graphs. Moreover, in this case, the normalizing
constant can be computed in closed form. We denote by $\pi_{U,\bolds{\alpha}}$
the normalized version of $\widetilde{\pi}_{U,\bolds{\alpha}}$ whenever
$z_G (U,\bolds{\alpha}) < \infty$. The following
lemma shows that the family $\pi_{U,\bolds{\alpha}}$ is a conjugate
family for Gaussian covariance graph models.

\begin{lemma}
Let $G = (V,E)$ be a decomposable graph, where vertices in $V$ are
ordered so that $V \in S_D$. Let ${\mathbf Y}_1, {\mathbf Y}_2, \ldots,
{\mathbf Y}_n$ be an i.i.d. sample from $\mathcal{N}_m ({\mathbf 0},
LDL^T)$, where $(L,D) \in {\bolds \Theta}_G$. Let $S = \frac{1}{n}
\sum_{i=1}^n {\mathbf Y}_i {\mathbf Y}_i^T$ denote the empirical
covariance matrix. If the prior distribution on $(L,D)$ is $\pi_{U,
\bolds{\alpha}}$, then the posterior distribution of $(L,D)$ is given
by $\pi_{\widetilde{U}, \widetilde{\bolds{\alpha}}}$, where
$\widetilde{U} = nS + U$ and $\widetilde{\bolds{\alpha}} = (n +
\alpha_1, n + \alpha_2, \ldots, n + \alpha_m)$.
\end{lemma}

\begin{pf}
The likelihood of the data is given by
\[
f({\mathbf y}_1, {\mathbf y}_2, \ldots, {\mathbf y}_n \mid L, D) =
\frac{1}{(\sqrt{2 \pi})^{nm}} e^{- ( \operatorname{tr} ( (LDL^T)^{-1} (nS) ) + n
\log |D| )/2}.
\]
Using $\pi_{U,\bolds{\alpha}}$ as a prior for $(L,D)$, the posterior
distribution of $(L,D)$ given the data $({\mathbf Y}_1, {\mathbf Y}_2,
\ldots, {\mathbf Y}_n)$ is
\begin{eqnarray*}
&&\pi_{U,\bolds{\alpha}} (L, D \mid {\mathbf Y}_1, {\mathbf Y}_2,
\ldots, {\mathbf Y}_n)\\
&&\qquad\propto e^{- ( \operatorname{tr}( (LDL^T)^{-1} (nS + U))
+ \sum_{i=1}^m (n + \alpha_i) \log D_{ii} )/2},\qquad\bolds{\theta}
\in {\bolds \Theta}_G.
\end{eqnarray*}
Hence, the posterior distribution belongs to the same family as the
prior, that is,
\[
\pi_{U,\bolds{\alpha}} ( \cdot \mid {\mathbf Y}_1, {\mathbf Y}_2,
\ldots, {\mathbf Y}_n) = \pi_{\widetilde{U}, \widetilde{\bolds{\alpha}}}
( \cdot ),
\]
where $\widetilde{U} = nS + U$ and $\widetilde{\bolds{\alpha}} = (n +
\alpha_1, n + \alpha_2, \ldots, n + \alpha_m)$.
\end{pf}

\begin{remark*}
If we assume that the observations have unknown mean
$\mu$, that is, ${\mathbf Y}_1, {\mathbf Y}_2, \ldots, {\mathbf Y}_n$
are i.i.d. $\mathcal{N} (\bolds{\mu}, \Sigma)$ with $\bolds{\mu} \in
\mathbb{R}^m,  \Sigma \in P_G$, then
\[
\widetilde{S} := \frac{1}{n} \sum_{i=1}^n ({\mathbf Y}_i - \bar{\mathbf Y})({\mathbf Y}_i - \bar{\mathbf Y})^T
\]
is the minimal sufficient statistic for $\Sigma$. Here, $n
\widetilde{S}$ has a Wishart distribution with parameter $\Sigma$ and
$n - 1$ degrees of freedom. Hence, if we assume a prior $\pi_{U,\bolds{\alpha}}$
for $(L,D)$, then the posterior distribution is given
by
\[
\pi_{U,\bolds{\alpha}} (\cdot \mid \widetilde{S}) =
\pi_{\widetilde{U}, \widetilde{\bolds{\alpha}}} (\cdot),
\]
where $\widetilde{U} = n \widetilde{S} + U$ and $\bolds{\alpha} =
(n-1+\alpha_1, n-1+\alpha_2, \ldots, n-1+\alpha_m)$.
\end{remark*}

\begin{remark*}
Note that, as with the distributions in
\cite{ltcmssmwdg,rajmasscar}, the functional form of the prior
distribution depends on the ordering of the vertices specified---but
this is not as restrictive as it first appears. In this sense, an
ordering is essentially another ``parameter'' to be specified and thus
can also be viewed as imposing extra information. We return to this
point in the examples section where we investigate the impact of
ordering on a real-world example (see Section \ref{examplessection}).
But, more importantly, given a perfect ordering of the vertices, any
rearrangement of the vertices within the residuals $R_j = C_j \backslash
H_{j-1}$ will still preserve the zeros between $\Sigma$ and $L,$ and
will thus be  sufficient for our purposes. In this sense, the
covariance Wishart distributions introduced in this paper do not
actually depend on a full ordering of the vertices. In fact, for the
class of decomposable graphs, any perfect ordering is sufficient, that
is, any ordering that is used in \cite{ltcmssmwdg} will also be
relevant for the covariance Wishart distributions defined above. In
this sense, these decomposable covariance Wishart distributions are
very flexible, especially since we are working in the curved
exponential family setting and are still able to use any ordering that
is appropriate for the \cite{ltcmssmwdg} distributions which address
the natural exponential family (NEF) concentration graph situation. The technical reason why any
perfect ordering will suffice is that any perfect ordering will
preserve the zeros between  $\Sigma$ and the matrix $L$ from its
Cholesky decomposition \cite{plsnpwrsmh,roverchldn}. Moreover, from an
applications perspective, since matrix operations are not invariant
with respect to ordering of the nodes, an ordering that facilitates
calculations is desirable. All that the ordering does is to relabel the
vertices, but the edge structure is completely and fully retained. To
further clarify what is meant, if one has a list of a genes/proteins
called ABLIM1, BCL6,  etc. and their names are replaced with the
numbers 1, 2, 3, etc., the problem can first be analyzed with the
integer labels and  one can then  go back to the original labels after
the analysis is done. So, in many applications, the ordering is not a
real restriction.
\end{remark*}

\subsection{Induced prior on $P_G$ and $Q_G$} \label{inducedprr}

The prior $\pi_{U, \bolds{\alpha}}$ on ${\bolds \Theta}_G$ (the
modified Cholesky space) induces a prior on $P_G$ (the covariance
matrix space) and $Q_G$. We provide an expression for the induced
priors on these spaces in order to compare our Wishart distributions
with other classes of distributions. Note that since the vertices have
been ordered so that $V \in S_D$, the transformation
\[
\psi\dvtx {\bolds \Theta}_G \rightarrow M_m^+
\]
defined by
\[
\psi (L,D) = LDL^T =: \Sigma
\]
is a bijection from ${\bolds \Theta}_G$ to $P_G$. The lemma below
provides the required Jacobians for deriving the induced priors on
$P_G$ and $Q_G$. The reader is referred to Section \ref{dmpgrphpty} for
notation on decomposable graphs. Note that if $x$ is a matrix, then
$|x|$ denotes its determinant, while if $C$ is a set, then $|C|$
denotes its cardinality.
\begin{lemma}[(Jacobians of transformations)]\label{jcbnttpqtp}

\begin{enumerate}
\item The Jacobian of the transformation $\psi\dvtx (L,D) \to \Sigma$ from
${\bolds \Theta}_G$ to $P_G$ is
\[
\prod_{i=1}^m D_{jj} (\Sigma)^{-n_j}.
\]
Here, $D_{jj} (\Sigma)$ denotes that $D_{jj}$ is a function of
$\Sigma$, and $n_j := |\{ i\dvtx (i,j) \in E, i > j \}|$ for $j = 1,2,
\ldots, m$.
\item The absolute value of the Jacobian of the bijection
$\zeta\dvtx x \rightarrow \hat{x}^{-1}$ from $Q_G$ to $P_G$ is
\[
\prod_{C \in \mathcal{C}} |x_C|^{- |C| - 1} \prod_{S \in \mathcal{S}}
|x_S|^{(|S| + 1) \nu (S)}.
\]
\end{enumerate}
\end{lemma}

\begin{pf}
The first part is a direct consequence of a result in
\cite{roveratond} and the proof of the second part can be found in
\cite{roverchldn}.
\end{pf}

These Jacobians allow us to compute the induced priors on $P_G$
and $Q_G$. The induced prior corresponding to
$\widetilde{\pi}_{U,\alpha}$ on $P_G$ is given by
\begin{equation} \label{distonpg}
\widetilde{\pi}^{P_G}_{U,\bolds{\alpha}} (\Sigma) \propto
e^{-(\operatorname{tr}(\Sigma^{-1} U) + \sum_{i=1}^m (2 n_i + \alpha_i) \log
D_{ii} (\Sigma) )/2},\qquad\Sigma \in P_G.
\end{equation}
We first note that the traditional inverse Wishart distribution (see
\cite{muirheadwd}) with parameters $U$ and $n$ is a special case of
(\ref{distonpg}) when $G$ is the complete graph and $\alpha_i = n - 2m
+ 2i,  \forall  1 \leq i \leq m$. We also note that the
$\mathcal{G}$-inverse Wishart priors introduced in \cite{ghrmnslvgs}
have a one-dimensional shape parameter $\delta$ and are a very special
case of our richer class $\widetilde{\pi}^{P_G}_{U,\bolds{\alpha}}$.
The single shape parameter $\delta$ is given by the relationship
$\alpha_i + 2 n_i = \delta + 2m,  1 \leq i \leq m$.\footnote{There is an
interesting parallel here that becomes apparent from our derivations
above. In the concentration graph setting, the single shape parameter
hyper-inverse Wishart (HIW) prior of Dawid and Lauritzen
\cite{dawidlrtzn} is a special case of  the multiple shape parameter
class of priors introduced by Letac and Massam \cite{ltcmssmwdg}, in
the sense that $\alpha_i = - \frac{1}{2} (\delta + c_i - 1)$ (see
\cite{rajmasscar} for notation). In a similar spirit, we discover that
the single shape parameter class of priors in \cite{ghrmnslvgs} is a
special case of the multiple shape parameter class of priors
$\widetilde{\pi}_{U,\bolds{\alpha}}$ introduced in this paper, in the
sense that $\alpha_i = \delta - 2 n_i + 2m$.}

We now proceed to derive the induced prior on $Q_G$. Let $x =
\kappa (\Sigma^{-1})$ denote the image of $\Sigma$ in $Q_G$ and let
$\hat{x}$ denote $\Sigma^{-1}$  (see \cite{ltcmssmwdg,rajmasscar} for
more details). Using the second part of Lemma \ref{jcbnttpqtp}, the
induced prior corresponding to $\widetilde{\pi}_{U,\bolds{\alpha}}$
on $Q_G$ is given by
\begin{eqnarray*}
\widetilde{\pi}^{Q_G}_{U,\bolds{\alpha}} (x)
&\propto& e^{-(\operatorname{tr}(\hat{x} U) + \sum_{i=1}^m (2 n_i + \alpha_i) \log D_{ii} ((\widehat{x})^{-1}))/2}\\
& & {}\times \frac{\prod_{S \in \mathcal{S}} |x_S|^{(|S|+1) \nu (S)}}{\prod_{C \in \mathcal{C}} |x_C|^{|C|+1}},\qquad x \in Q_G.
\end{eqnarray*}
%
\subsection{Comparison with the Letac--Massam priors} \label{ltcmssmcpn}

We now carefully compare our class of priors to those proposed in Letac
and Massam \cite{ltcmssmwdg}. In \cite{ltcmssmwdg}, the authors
construct two classes of distributions, named $W_{P_G}$ and $W_{Q_G}$,
on the spaces $P_G$ and $Q_G$, respectively, for $G$ decomposable (see
\cite{ltcmssmwdg}, Section 3.1). These distributions are
generalizations of the Wishart distribution on these convex cones and
have been found to be very useful for high-dimensional Bayesian
inference, as illustrated in \cite{rajmasscar}. These priors lead to
corresponding classes of inverse Wishart distributions $\mathit{IW}_{P_G}$ (on
$Q_G$) and $\mathit{IW}_{Q_G}$ (on $P_G$), that is, $U \sim \mathit{IW}_{P_G}$ whenever
$\hat{U}^{-1} \sim W_{P_G}$, and $V \sim \mathit{IW}_{Q_G}$ whenever $\kappa
(V^{-1}) \sim W_{Q_G}$. In \cite{ltcmssmwdg}, it is shown that the
family of distributions $\mathit{IW}_{P_G}$ yields a family of conjugate priors
in the Gaussian concentration graph setting, that is, when $\Sigma \in
Q_G$.

As the $\mathit{IW}_{Q_G}$ priors of \cite{ltcmssmwdg} are defined on
the space $P_G$, in principle, they can potentially serve as
priors\footnote{The use of this class of nonconjugate priors for
Bayesian inference in covariance graph models was already explored in
\cite{lmwdgvpcvg}.} in the covariance graph setting since the parameter
of interest $\Sigma$ lives in $P_G$. Let us examine this class more
carefully, first with a view to understanding their use in the
covariance graph setting and second to compare them to our priors.
Following the notation for decomposable graphs in Section
\ref{dmpgrphpty} and in \cite{ltcmssmwdg}, the density of the
$\mathit{IW}_{Q_G}$ distribution is given by
\[
{\mathit{IW}}_{Q_G}^{U,\bolds{\alpha}, \bolds{\beta}} (\Sigma) \propto
e^{\operatorname{tr} ( \Sigma^{-1} U )/2} \frac{\prod_{C \in \mathcal{C}} |
(\Sigma^{-1})_C |^{\alpha (C)+ (c+1)/2}} {\prod_{S \in
\mathcal{S}} | (\Sigma^{-1})_S |^{\nu (S)(\beta (S) + (s+1)/2)}},
\qquad\Sigma \in P_G,
\]
where $U \in P_G$, and $\alpha (C)$, $C \in \mathcal{C}$ and $\beta
(S)$, $S \in \mathcal{S}$ are real numbers. The posterior density of
$\Sigma$ under this prior is given by
\begin{eqnarray*}
\pi^{\mathrm{IW}}_{U,\bolds{\alpha},\bolds{\beta}} (\Sigma \mid {\mathbf
Y}_1, {\mathbf Y}_2, \ldots, {\mathbf Y}_n) &\propto&
e^{-\operatorname{tr}
( \Sigma^{-1} (U + nS) )/2}\\
&&{}\times\frac{\prod_{C \in \mathcal{C}} |
(\Sigma^{-1})_C |^{\alpha (C)+ (c+1)/2 + n/2}} {\prod_{S
\in \mathcal{S}} | (\Sigma^{-1})_S |^{\nu (S)(\beta (S) +
(s+1)/2) + n \nu (S)/2}}.
\end{eqnarray*}
However, $U + nS$ may not, in general, be in $P_G$, which is a crucial
assumption in the analysis in \cite{ltcmssmwdg}. Hence, the conjugacy
breaks down.

We now investigate similarities and differences between our
class of priors and the $\mathit{IW}_{Q_G}$ class. Since the $\mathit{IW}_{Q_G}^{U,
\bolds{\alpha}, \bolds{\beta}}$ density\vspace*{1pt} is defined only for $U \in
P_G$, a pertinent question is whether our class of priors has the same
functional form when $U \in P_G$. We discover that this is not the case
and demonstrate this through an example. Consider the $4$-chain $A_4$.
One can easily verify that the terms $e^{-\operatorname{tr}( \Sigma^{-1}
U )/2}$ are identical in both priors. We now show that the remaining
terms are not identical. If $\Sigma = LDL^T$ is the modified Cholesky
decomposition of $\Sigma$, then, for this particular graph with $C_1 =
\{1,2\}$, $C_2 = \{2,3\}$, $C_3 = \{3,4\}$ and $ S_2 = \{3\}$, $S_3 =
\{4\}$, the expression that is not in the exponential term for the
$\mathit{IW}_{Q_G}$ density is of the form
\begin{eqnarray*}
\frac{\prod_{i=1}^3 | (\Sigma^{-1})_{C_i} |^{\alpha_i}} {\prod_{i=2}^3
| (\Sigma^{-1})_{S_i} |^{\beta_i}} &=& \biggl( \frac{1}{D_{11}} \biggr)^{\alpha_1}
\biggl( \frac{1}{D_{22}} +
\frac{L_{32}^2}{D_{33}} + \frac{L_{32}^2 L_{43}^2}{D_{44}} \biggr)^{\alpha_1 - \beta_1}\\
& &{} \times \biggl( \frac{1}{D_{22}}\biggr)^{\alpha_2} \biggl( \frac{1}{D_{33}} +
\frac{L_{43}^2}{D_{44}} \biggr)^{\alpha_2 - \beta_2} \biggl( \frac{1}{D_{33}
D_{44}} \biggr)^{\alpha_3}.
\end{eqnarray*}
This expression is clearly different from the term, other than the
exponent $e^{-\operatorname{tr}( \Sigma^{-1} U )/2}$ in $\pi^{P_G}_{U,
\bolds{\alpha}}$, which is a product of different powers of $D_{ii},
 i = 1, 2, 3, 4$.

However, an interesting property emerges when $G$ is
homogeneous. Note that, in this case, for any clique $C$ and any
separator $S$,
\[
| (\Sigma^{-1})_C | = \prod_{i \in C} \frac{1}{D_{ii}},   |
(\Sigma^{-1})_S | = \prod_{i \in S} \frac{1}{D_{ii}}.
\]
Hence, when $G$ is homogeneous, the class $\mathit{IW}_{Q_G}$ is contained in
the class $\pi^{P_G}$. The containment is strict because $U$ need not
be in $P_G$ for  our class $\pi^{P_G}$. Also, in $\mathit{IW}_{Q_G}$, the
exponent of $D_{ii}$ and $D_{jj}$ is the same if $i R j$, that is, the
shape parameter, is shared for vertices in the same equivalence class,
as defined by the relation $R$. We, however, note that the difference
in the number of shape parameters is not a major difference, due to the
result of Consonni and Veronese \cite{cnsnvrnscr}, together with fact
that for the $W_{Q_G}$ (and, correspondingly, for the $\mathit{IW}_{Q_G}$), each
one of the blocks $x_{[i]\cdot}$ has a Wishart distribution (see Theorem 4.5 of
\cite{ltcmssmwdg}).

We therefore note that in the restrictive case when $G$ is homogeneous
and when $U \in P_G$, the two classes of distributions $\pi^{P_G}$ and
$\mathit{IW}_{Q_G}$ have the same functional form. The fact that we do not
restrict $U \in P_G$ is an important difference since, even in the
homogeneous case, they yield  a larger class of distributions on the
homogeneous cone $P_G$ compared to those in Andersson and Wojnar
\cite{adscnvcmcs}, resulting in nonsuperficial consequences for
inference in covariance graph models.\footnote{We note, however, that
the distributions in \cite{adscnvcmcs} are quite general since the
authors  consider other homogeneous cones and not just $P_G$.}

\section{Sampling from the posterior distribution} \label{cndtlcmptn}

In this section, we study the properties of our family of distributions
and thereby provide a method that allows us to generate samples from
the posterior distribution corresponding to the priors defined in
Section \ref{priorcnsrt}. In particular, we prove that $\bolds{\theta}
= (L,D) \in {\bolds \Theta}_G$ can be partitioned into blocks so that
the conditional distribution of each block given the others is a
standard distribution in statistics and hence easy to sample from. We
can therefore generate samples from the posterior distribution by using
the block Gibbs sampling algorithm.

\subsection{Distributional properties and the block Gibbs sampler} \label{blkgssimul}

Let us introduce some notation before deriving the required conditional
distributions. Let $G = (V,E)$ be a decomposable graph such that $V \in
S_D$. For a lower-triangular matrix $L$ with diagonal entries equal to
$1$,
\begin{eqnarray*}
L_{u \cdot} &:=& u\mbox{th row of } L,\qquad u = 1, 2, \ldots, m,\\
L_{\cdot v} &:=& v\mbox{th column of } L,\qquad v = 1, 2, \ldots, m,\\
L^G_{\cdot v} &:=& (L_{uv})_{u>v, (u,v) \in E},\qquad v = 1, 2, \ldots,
m-1.
\end{eqnarray*}
So, $L^G_{\cdot v}$ is the $v$th column of $L$ without the components
which are specified to be zero under the model $\mathcal{G}$ (and
without the $v$th diagonal entry, which is $1$). In terms of this
notation, the parameter space can be represented as
\begin{eqnarray} \label{coldcmpmcl}
&&{\bolds \Theta}_G = \{ ( L^G_{\cdot 1}, L^G_{\cdot 2}, L^G_{\cdot 3},
\ldots, L^G_{\cdot m-1}, D )\dvtx\nonumber\\ [-8pt]\\ [-8pt]
&&\hphantom{{\bolds \Theta}_G = \{}L_{ij} \in \mathbb{R},   \forall  1 \leq
j < i \leq m, (i,j) \in E,  D_{ii}> 0,  \forall  1 \leq i \leq m
\}.\nonumber
\end{eqnarray}
Suppose that $\bolds{\theta}\sim \pi_{U, \alpha}$ for some positive
definite $U$ and $\bolds{\alpha} \in \mathbb{R}^m$ with nonnegative
entries. The posterior distribution is then $\pi_{\widetilde{U},
\widetilde{\bolds{\alpha}}}$, where $\widetilde{U} = nS + U,
\widetilde{\bolds{\alpha}} = (n + \alpha_1, n + \alpha_2, \ldots, n +
\alpha_m)$. In the following proposition, we derive the distributional
properties which provide the essential ingredients for our block Gibbs
sampling procedure.

\begin{theorem} \label{cndtldstbn}
Using the notation above, the conditional distributions of each
component of $\bolds{\theta}$ [as in (\ref{coldcmpmcl})] given the
other components and the data ${\mathbf Y}_1, {\mathbf Y}_2, \ldots,
{\mathbf Y}_n$ are as follows:
\begin{enumerate}
\item
\[
L^G_{\cdot v} \mid (L \backslash L^G_{\cdot v}, D, {\mathbf Y}_1,
{\mathbf Y}_2, \ldots, {\mathbf Y}_n) \sim \mathcal{N} (
\bolds{\mu}^{v,G}, M^{v,G})\qquad\forall  v = 1, 2, \ldots, m-1,
\]
where
\begin{eqnarray*}
\mu^{v,G}_u &:=& \mu^v_u + \sum_{u'>v: (u',v) \in E} \sum_{w>v: (w,v)
\notin E \atop{\mathit{or}\  w < v, L^{-1}_{vw} = 0}} M^{v,G}_{uu'} (
L^{-1} \widetilde{U}
(L^T)^{-1} )_{vv}\\
&&\hphantom{\mu^v_u + \sum_{u'>v: (u',v) \in E} \sum_{w>v: (w,v)
\notin E \atop{\mathrm{or}\ w < v, L^{-1}_{vw} = 0}}}{}\times(LDL^T)^{-1}_{u'w} \mu^v_w\\
&&\hspace*{200pt}\forall u > v, (u,v) \in E,\\
\mu^v_u &:=& \frac{(L^{-1} \widetilde{U})_{vu}}{( L^{-1} \widetilde{U}
(L^T)^{-1})_{vv}}\qquad\forall u \mbox{ such that } L^{-1}_{vu} = 0,\\
(M^{v,G})^{-1}_{uu'} &:=& ( L^{-1} \widetilde{U} (L^T)^{-1} )_{vv}
(LDL^T)^{-1}_{uu'}\qquad\forall u,u' > v, (u,v),(u',v) \in E;
\end{eqnarray*}
\item
\[
D_{ii} \mid L, {\mathbf Y}_1, {\mathbf Y}_2, \ldots, {\mathbf Y}_m \sim
\operatorname{IG}\biggl( \frac{\widetilde{\alpha}_i}{2} - 1, \frac{( L^{-1} \widetilde{U}
(L^T)^{-1} )_{ii}}{2} \biggr),
\]
independently for $i = 1, 2, \ldots, m$, where ``$\operatorname{IG}$'' represents the
inverse-gamma distribution.
\end{enumerate}
\end{theorem}

\begin{remark*}
The notation $w: L^{-1}_{vw} = 0$ in the definition of
$\bolds{\mu}^{v,G}$ above means indices $w$ for which $L^{-1}_{vw}$ is
$0$ as a function of entries of $L$.
\end{remark*}

Deriving the required conditional distributions in Theorem
\ref{cndtldstbn} entails careful analysis. We first state two lemmas
which are essential for deriving these distributions.

\begin{lemma} \label{derivative}
Let $u > v$, $(u,v) \in E$. Then,
\[
\frac{\partial L^{-1}_{ij}}{\partial L_{uv}} = - L^{-1}_{iu}
L^{-1}_{vj}\qquad\forall 1 \leq j < i \leq m.
\]
\end{lemma}

\begin{pf}
The proof is straightforward and is therefore omitted for
brevity.
\end{pf}

Recall from Proposition \ref{inverselrt} that $L^{-1}_{ij}$
functionally depends on $L_{uv}$ only if $i \geq u > v \geq j$. We use
this observation repeatedly in our arguments. For a given $v$, to prove
conditional multivariate normality of the  conditional distribution of
$L^G_{\cdot v}$ given the others, we shall demonstrate that if we treat
$D$ and the other columns of $L$ as constants, then $\operatorname{tr}((LDL^T)^{-1}
\widetilde{U})$ is a quadratic form in the entries of $L^G_{\cdot v}$ .

\begin{lemma} \label{hesscalcpr}
Let $u,u' > v$, $(u,v),(u',v) \in E$. Then,
\[
\frac{\partial^2}{\partial L_{uv}\,\partial L_{u'v}} \operatorname{tr}((LDL^T)^{-1}
\widetilde{U}) = 2 ( L^{-1} \widetilde{U} (L^T)^{-1} )_{vv}
(LDL^T)^{-1}_{uu'},
\]
which is functionally independent of the elements of $L^G_{\cdot v}$.
\end{lemma}

\begin{pf} First, note that,
\begin{eqnarray*}
& & \frac{\partial}{\partial L_{uv}} \operatorname{tr}((LDL^T)^{-1} \widetilde{U})\\
&&\qquad= \frac{\partial}{\partial L_{uv}} \Biggl(\sum_{i=1}^m \sum_{j=1}^m
\sum_{k=1}^m\frac{L^{-1}_{ki} L^{-1}_{kj}}{D_{kk}} \widetilde{U}_{ij} \Biggr)\\
&&\qquad= - \sum_{i=1}^m \sum_{j=1}^m \sum_{k=1}^m \biggl(\frac{L^{-1}_{ku}
L^{-1}_{vi}L^{-1}_{kj} + L^{-1}_{ki} L^{-1}_{ku} L^{-1}_{vj}}{D_{kk}}\biggr) \widetilde{U}_{ij}\qquad(\mbox{by Lemma } \ref{derivative})\\
&&\qquad= -2 \sum_{i=1}^m \sum_{j=1}^m \sum_{k=1}^m \frac{L^{-1}_{ku}
L^{-1}_{vi} L^{-1}_{kj}} {D_{kk}} \widetilde{U}_{ij}.
\end{eqnarray*}
Note that $L^{-1}$ is a lower-triangular matrix. Hence,
\begin{eqnarray*}
&& \frac{\partial^2}{\partial L_{uv}\,\partial L_{u'v}} \operatorname{tr}((LDL^T)^{-1} \widetilde{U})\\
&&\qquad = -2 \frac{\partial}{\partial L_{u'v}} \Biggl( \sum_{i=1}^m \sum_{j=1}^m\sum_{k=1}^m
\frac{L^{-1}_{ku} L^{-1}_{vi} L^{-1}_{kj}}{D_{kk}} \widetilde{U}_{ij}\Biggr)\\
&&\qquad = 2 \sum_{i=1}^m \sum_{j=1}^m \sum_{k=1}^m \frac{L^{-1}_{ku}
L^{-1}_{vi} L^{-1}_{ku'}L^{-1}_{vj}}{D_{kk}} \widetilde{U}_{ij}\\
&&\qquad = 2\Biggl(\sum_{i=1}^m \sum_{j=1}^m L^{-1}_{vi} \widetilde{U}_{ij}L^{-1}_{vj}\Biggr)
\Biggl(\sum_{k=1}^m \frac{L^{-1}_{ku} L^{-1}_{ku'}}{D_{kk}}\Biggr)\\
&&\qquad = 2 ( L^{-1} \widetilde{U} (L^T)^{-1} )_{vv} (LDL^T)^{-1}_{uu'}.
\end{eqnarray*}
The second equality above follows by noting that by Proposition
\ref{inverselrt}, $L^{-1}_{vi}$ is functionally independent of
$L^G_{\cdot v}$ for all $1 \leq i \leq m$ and $L^{-1}_{ku}$ is
functionally independent of $L^G_{\cdot v}$ for all $1 \leq k \leq m$
and $u > v$, and then applying Lemma \ref{derivative}. Using this
functional independence argument above once more, we thereby conclude
that $2 ( L^{-1} \widetilde{U} (L^T)^{-1} )_{vv} (LDL^T)^{-1}_{uu'}$ is
independent of $L^G_{\cdot v}$.
\end{pf}

\begin{pf*}{Proof of Theorem \ref{cndtldstbn}}
An immediate consequence of
Lemma \ref{hesscalcpr} and the\vspace*{1pt} preceding remark is that we can write
$\operatorname{tr}((LDL^T)^{-1} \widetilde{U})$ as follows:
\begin{eqnarray*}
&& \operatorname{tr}((LDL^T)^{-1} \widetilde{U})\\
&&\qquad= \sum_{u>v, (u,v) \in E} \sum_{u'>v, (u',v) \in E} ( ( L^{-1}
\widetilde{U} (L^T)^{-1} )_{vv} (LDL^T)^{-1}_{uu'} )\\
&&\qquad\hphantom{= \sum_{u>v, (u,v) \in E} \sum_{u'>v, (u',v) \in E}}{}\times(L_{uv} - b_u)
(L_{u'v} - b_{u'}) + C,
\end{eqnarray*}
where ${\mathbf b} = (b_u)_{u>v, (u,v) \in E}$ and $C$ are independent
of $L^G_{\cdot v}$. In order to evaluate $(b_u)_{u>v, (u,v) \in E}$,
note that the term in $\frac{\partial}{\partial L_{uv}} \operatorname{tr}((LDL^T)^{-1}
\widetilde{U})$ which is independent of $L^G_{\cdot v}$ is given by
\begin{equation} \label{exprconst1}
-2 \sum_{u'>v, (u',v) \in E} ( ( L^{-1} \widetilde{U} (L^T)^{-1} )_{vv}
(LDL^T)^{-1}_{uu'} ) b_{u'}
\end{equation}
for every $u>v$, $(u,v) \in E$. However, from the proof of Lemma
\ref{hesscalcpr}, we alternatively know that
\[
\frac{\partial}{\partial L_{uv}} \operatorname{tr}((LDL^T)^{-1} \widetilde{U}) = -2
\sum_{i=1}^m \sum_{j=1}^m \sum_{k=1}^m \frac{L^{-1}_{ku} L^{-1}_{vi}
L^{-1}_{kj}}{D_{kk}} \widetilde{U}_{ij}.
\]
Note that by Lemma \ref{derivative}, $L^{-1}_{ku} L^{-1}_{kj}$ is
functionally dependent on $L^G_{\cdot v}$ if and only if $L^{-1}_{ku}
\neq 0$ and $L^{-1}_{vj} \neq 0$ (as a function of $L$). Hence, the
term in $\frac{\partial}{\partial L_{uv}} \operatorname{tr}((LDL^T)^{-1}
\widetilde{U})$ which is independent of $L^G_{\cdot v}$ is given by
%
\begin{eqnarray}\label{exprconst2}
&& -2 \sum_{i=1}^m \sum_{j=1}^m \sum_{k=1}^m \frac{L^{-1}_{ku}
L^{-1}_{vi} L^{-1}_{kj}}{D_{kk}}
\widetilde{U}_{ij} 1_{\{L^{-1}_{vj} = 0\ \mathrm{or}\ L^{-1}_{ku} = 0\}} \nonumber\\
&&\qquad= -2 \sum_{j: L^{-1}_{vj}=0}\Biggl( \sum_{i=1}^m L^{-1}_{vi}
\widetilde{U}_{ij}\Biggr)
\Biggl( \sum_{k=1}^m \frac{L^{-1}_{ku} L^{-1}_{kj}}{D_{kk}} \Biggr)\\
&&\qquad= -2 \sum_{j: L^{-1}_{vj}=0} (L^{-1} \widetilde{U} )_{vj}
(LDL^T)^{-1}_{uj}\nonumber\\
&&\qquad= -2 \sum_{j: L^{-1}_{vj}=0} \frac{( L^{-1} \widetilde{U} )_{vj}}{(
L^{-1} \widetilde{U} (L^T)^{-1} )_{vv}} ( L^{-1} \widetilde{U}
(L^T)^{-1} )_{vv} (LDL^T)^{-1}_{uj}.\nonumber
\end{eqnarray}
Now, observe the following facts:
\begin{enumerate}
\item the expressions in (\ref{exprconst1}) and (\ref{exprconst2})
should be the same for every $u > v$, $(u,v) \in E$;
\item if $A = \left({{A_1\enskip  A_2}\atop {A_2^T\enskip A_3}}\right),
\bolds{\xi} = \left({{\bolds{\xi}_1}\atop {\bolds{\xi}_2}}\right)$ and $\bolds{\eta}$ are such that
\[
A_1 \bolds{\xi}_1 + A_2 \bolds{\xi}_2 = A_1 \bolds{\eta},
\]
then
\[
\bolds{\eta} = \bolds{\xi}_1 + A_1^{-1} A_2 \bolds{\xi}_2.
\]
\end{enumerate}
If we choose $A$, $\bolds{\xi}$ and $\bolds{\eta}$ as
\begin{eqnarray*}
A_{uu'} &:=& ( L^{-1} \widetilde{U} (L^T)^{-1} )_{vv}
(LDL^T)^{-1}_{uu'}\qquad\forall u,u' \mbox{ such
that } L^{-1}_{vu}, L^{-1}_{vu'} = 0,\\
\xi_u &:=& \frac{(L^{-1} \widetilde{U})_{vu}}{( L^{-1} \widetilde{U}
(L^T)^{-1} )_{vv}}\qquad\forall u
\mbox{ such that } L^{-1}_{vu} = 0,\\
\eta_u &:=& b_u \qquad\forall u > v, (u,v) \in E,
\end{eqnarray*}
then combining the observations above with (\ref{exprconst1}) and
(\ref{exprconst2}), we obtain
\begin{eqnarray*}
&&\operatorname{tr}((LDL^T)^{-1} \widetilde{U})\\
&&\qquad= \sum_{u>v, (u,v) \in E} \sum_{u'>v, (u',v) \in E} ( ( L^{-1}
\widetilde{U} (L^T)^{-1} )_{vv} (LDL^T)^{-1}_{uu'} )\\
&&\qquad\hphantom{= \sum_{u>v, (u,v) \in E} \sum_{u'>v, (u',v) \in E}}{}\times(L_{uv} -
\mu^{v,G}_u)(L_{u'v} - \mu^{v,G}_{u'}) + C.
\end{eqnarray*}
As defined earlier,
%
\begin{eqnarray}
\mu^v_u &=& \frac{( L^{-1} \widetilde{U} )_{vu}}{( L^{-1} \widetilde{U}
(L^T)^{-1})_{vv}} \qquad\forall u \mbox{ such that } L^{-1}_{vu} = 0,\nonumber\\
\mu^{v,G}_u &=& \mu^v_u + \sum_{u'>v, (u',v) \in E} \sum_{w>v, (w,v)
\notin E \atop{\ \mathrm{or}\  w < v, L^{-1}_{vw} = 0}}
M^{v,G}_{uu'} ( L^{-1} \widetilde{U}
(L^T)^{-1} )_{vv}(LDL^T)^{-1}_{u'w} \mu^v_w\nonumber\\
\eqntext{\forall u>v, (u,v) \in
E}
\end{eqnarray}
and $C$ is independent of $L^G_{\cdot v}$. It follows that under
$\pi_{\widetilde{U}, \widetilde{\alpha}}$, the conditional distribution
of $L^G_{\cdot v}$ given the other parameters and the data ${\mathbf
Y}_1, {\mathbf Y}_2, \ldots, {\mathbf Y}_n$ is $\mathcal{N} (
\bolds{\mu}^{v,G}, M^{v,G})$.

For deriving the conditional distribution of the entries of
$D$, we note that
\[
e^{- ( \operatorname{tr}((LDL^T)^{-1} \widetilde{U}) + \sum_{i=1}^m
\widetilde{\alpha}_i \log D_{ii} )/2} = \prod_{j=1}^m
\frac{1}{D_{jj}^{\widetilde{\alpha}_j/2}} e^{- ( L^{-1}
\widetilde{U} (L^T)^{-1} )_{jj}/(2 D_{jj})}.
\]
The above leads us to conclude that the conditional distribution of
$D_{jj}$ given the other parameters and the data ${\mathbf Y}_1,
{\mathbf Y}_2, \ldots, {\mathbf Y}_n$ is $\operatorname{IG} (
\frac{\widetilde{\alpha}_j}{2} - 1, \frac{( L^{-1} \widetilde{U}
(L^T)^{-1} )_{jj}}{2} )$, independently for every $j = 1, 2, \ldots,m$.
\end{pf*}

\subsection{Convergence of block Gibbs sampler} \label{gibbscnvgc}

The block Gibbs sampling procedure, based on the conditional
distributions derived above, gives rise to a Markov chain. It is
natural to ask whether this Markov chain converges to the desired
distribution $\pi_{\widetilde{U}, \widetilde{\bolds{\alpha}}}$.
Convergence properties are sometimes overlooked due to the theoretical
demands in establishing them. However, they yield theoretical
safeguards that the block Gibbs sampling algorithm can be used for
sampling from the posterior distribution.

We now prove that sufficient conditions for convergence of a
Gibbs sampling Markov chain to its stationary distribution (see
\cite{adscnvcmcs}, Theorem 6) are satisfied by the Markov chain
corresponding to our block Gibbs sampler. Let $\phi ({\mathbf x} \mid
\bolds{\mu}, \Sigma)$ denote the $\mathcal{N} (\bolds{\mu}, \Sigma)$
density evaluated at ${\mathbf x}$. Let $f_{\operatorname{IG}} (d \mid \alpha,
\lambda)$ denote the $\operatorname{IG} (\alpha, \lambda)$ density evaluated at $d$.
Let us fix $\psi, d_1, d_2 > 0$ arbitrarily. Let
\[
{\bolds \Theta}_{\psi, d_1, d_2} := \{ \bolds{\theta} = (L,D) \in
{\bolds \Theta}_G\dvtx   |L_{ij}| \leq \psi, d_1 \leq D_{ii} \leq d_2\
\forall  i > j, (i,j) \in E\}.
\]
We now formally prove the conditions which are sufficient for
establishing convergence.

\begin{prop} \label{adsgscndtn}
There exists some $ \delta > 0$ such that,  uniformly for all $
\bolds{\theta }= (L,D) \in {\bolds \Theta}_{\psi, d_1, d_2}$,
\begin{eqnarray*}
 \phi ( L^G_{\cdot, v} \mid \bolds{\mu}^{v,G}, M^{v,G} ) &>& \delta\qquad
 \forall  v = 1,2, \ldots,m-1,\\
f_{\operatorname{IG}} \biggl( D_{ii} \big\vert \frac{\widetilde{\alpha}_i}{2} - 1, \frac{(
L^{-1} \widetilde{U} (L^T)^{-1} )_{ii}}{2} \biggr) &>& \delta  \qquad\forall  i =
1,2, \ldots, m.
\end{eqnarray*}
\end{prop}

\begin{pf}First, by Proposition \ref{inverselrt}, all entries of
$L^{-1}$ are polynomials in the entries of $L$. Since ${\bolds
\Theta}_{\psi, d_1, d_2}$ is bounded and closed, there exists $\psi_1 >
0$ such that
\[
(L,D) \in {\bolds \Theta}_{\psi, d_1, d_2} \Rightarrow | L^{-1}_{uv} |
\leq \psi_1\qquad\forall u > v, (u,v) \in E.
\]
Using the above, there exists a constant $\psi_2 > 0$ such that if
$(L,D) \in {\bolds \Theta}_{\psi, d_1, d_2}$, then
\begin{equation} \label{estimates1}
\hspace*{24pt}| ( L^{-1} \widetilde{U} )_{vu} | \leq \psi_2,\qquad| (LDL^T)^{-1}_{uu'} |
\leq \psi_2,\qquad| ( L^{-1} \widetilde{U} (L^T)^{-1} )_{vv} | \leq \psi_2
\end{equation}
for every $1 \leq v,u,u' \leq m$. Second, since $L^{-1}_{vv} = 1$ for
all $ 1 \leq v \leq m$ and $\widetilde{U}$ is positive definite, it
follows that there exists a constant $\psi_3 > 0$ such that if $(L,D)
\in {\bolds \Theta}_{\psi, d_1, d_2}$, then
\begin{equation} \label{estimates2}
\psi_3 \leq ( L^{-1} \widetilde{U} (L^T)^{-1} )_{vv}
\end{equation}
for every $1 \leq v \leq m$.

Let $(L,D) \in {\bolds \Theta}_{\psi, d_1, d_2}$. Note that
\begin{eqnarray*}
& & ( L^G_{\cdot v} - \bolds{\mu}^{v,G} )^T ( M^{v,G} )^{-1} (
L^G_{\cdot v} - \mu^{v,G} )\\
&&\qquad= ( L^G_{\cdot v} )^T ( M^{v,G} )^{-1} L^G_{\cdot v} - 2 ( L^G_{\cdot
v} )^T ( M^{v,G} )^{-1} \mu^{v,G} + ( \mu^{v,G} )^T ( M^{v,G} )^{-1}
\mu^{v,G}.
\end{eqnarray*}
Observe that if $ \bolds{\zeta} = \left({{\bolds{\zeta}_1}\atop {\bolds{\zeta}_2}}\right) \in \mathbb{R}^m$ and $\Sigma = \left({{\Sigma_{11}\enskip \Sigma_{12}}
\atop {\Sigma_{21} \enskip \Sigma_{22}}}\right)$ is a positive definite
matrix, then
\[
( \bolds{\zeta}_1 + \Sigma_{11}^{-1} \Sigma_{12} \bolds{\zeta}_2 )^T
\Sigma_{11} ( \bolds{\zeta}_1 + \Sigma_{11}^{-1} \Sigma_{12}
\bolds{\zeta}_2 ) \leq \zeta^T \Sigma \zeta \qquad \forall  \zeta \in \mathbb{R}^m.
\]
If we choose $\bolds{\zeta}$ and $\Sigma$ as
\begin{eqnarray*}
\Sigma_{uu'} &:=& ( L^{-1} \widetilde{U} (L^T)^{-1} )_{vv}
(LDL^T)^{-1}_{uu'}\qquad \forall u,u'
\mbox{ such that } L^{-1}_{vu}, L^{-1}_{vu'} = 0,\\
\zeta_u &=& \mu^v_u \qquad\forall u \mbox{ such that } L^{-1}_{vu} = 0,
\end{eqnarray*}
then combining the observation above and the definition of $\bolds{\mu}^{v,G}$, we get that
\[
( \bolds{\mu}^{v,G} )^T ( M^{v,G} )^{-1} \bolds{\mu}^{v,G} \leq
\sum_{u: L^{-1}_{vu} = 0} \sum_{u': L^{-1}_{vu'} = 0} \mu^v_u ( L^{-1}
\widetilde{U} (L^T)^{-1} )_{vv} (LDL^T)^{-1}_{uu'} \mu^v_{u'}.
\]
From the definitions in Theorem \ref{cndtldstbn}, we also have
\[
( ( M^{v,G} )^{-1} \bolds{\mu}^{v,G} )_u = \sum_{j: L^{-1}_{vj} = 0} (
L^{-1} \widetilde{U} )_{vj} (LDL^T)^{-1}_{uj}\qquad\forall u > v, (u,v) \in
E.
\]
It follows by (\ref{estimates1}) that for $(L,D) \in {\bolds
\Theta}_{\psi, d_1, d_2}$, there exists $\psi_4 > 0$ such that
\begin{equation} \label{estimates3}
( L^G_{\cdot v} - \bolds{\mu}^{v,G} )^T ( M^{v,G} )^{-1} ( L^G_{\cdot
v} - \bolds{\mu}^{v,G} ) \leq \psi_4
\end{equation}
for every $v = 1,2, \ldots, m-1$. Also, by the definition of $M^{v,G}$,
it follows that for $(L,D) \in {\bolds \Theta}_{\psi, d_1, d_2},  0 < |
M^{v,G} | < \infty$ and $| M^{v,G} |$ is a continuous function of
$(L,D)$. Recall that for a matrix $A,  |A|$ denotes the determinant of
$A$. Since ${\bolds \Theta}_{\psi, d_1, d_2}$ is a bounded and closed
set, both the maximum and minimum of the function $| M^{v,G} |$ are
attained in ${\bolds \Theta}_{\psi, d_1, d_2}$. It follows that for
$(L,D) \in {\bolds \Theta}_{\psi, d_1, d_2}$, there exist $0 < \kappa_1
< \kappa_2$ such that
\begin{equation} \label{estimates4}
\kappa_1 < | M^{v,G} | < \kappa_2
\end{equation}
for every $v = 1,2, \ldots, m - 1$. It follows by (\ref{estimates1}),
(\ref{estimates2}), (\ref{estimates3}) and (\ref{estimates4}) that for
$(L,D) \in {\bolds \Theta}_{\psi, d_1, d_2}$, there exists $\delta_1 >
0$ such that
\[
\phi ( L^G_{\cdot v} \mid \bolds{\mu}^{v,G}, M^{v,G} ) > \delta_1
\qquad\forall  v = 1,2, \ldots, m-1.
\]

Note, furthermore, that if $(L,D) \in {\bolds \Theta}_{\psi, d_1,
d_2}$, then, from (\ref{estimates1}) and (\ref{estimates2}),
\[
\psi_3 \leq ( L^{-1} \widetilde{U} (L^T)^{-1} )_{ii} \leq \psi_2
\]
for every $1 \leq i \leq m$. Hence, there exists $\delta_2 > 0$ such
that
\[
f_{\operatorname{IG}} \biggl( D_{ii} \big\vert \frac{\widetilde{\alpha}_i}{2} - 1, \frac{(
L^{-1} \widetilde{U} (L^T)^{-1})_{ii}}{2}\biggr ) > \delta_2\qquad\forall  i =
1,2, \ldots, m.
\]
Let $\delta = \min (\delta_1, \delta_2)$. Hence, for $\bolds{\theta} =
(L,D) \in {\bolds \Theta}_{\psi, d_1, d_2}$,
\begin{eqnarray*}
\phi ( L^G_{\cdot v} \mid \bolds{\mu}^{v,G}, M^{v,G} ) &>& \delta
\qquad\forall  v= 1,2, \ldots, m-1,\\
f_{\operatorname{IG}} \biggl( D_{ii} \big\vert \frac{\widetilde{\alpha}_i}{2} - 1, \frac{(
L^{-1} \widetilde{U} (L^T)^{-1} )_{ii}}{2} \biggr) &>& \delta \qquad\forall  i =
1,2, \ldots, m.
\end{eqnarray*}
\upqed\hspace*{-40pt}\end{pf}

Recall that $n_v = |\{u\dvtx   u > v, (u,v) \in E\}|$. Note that the
measures corresponding to $\mathcal{N} ( \bolds{\mu}^{v,G}, M^{v,G} )$
and $\mathcal{N} ( {\mathbf 0}, I_{n_v} )$ are mutually absolutely
continuous and the corresponding densities with respect to Lebesgue
measure are positive everywhere on $\mathbb{R}^{n_v}$ for all $v = 1,2,
\ldots, m-1$. In addition, the measures corresponding to $\operatorname{IG} (
\frac{\widetilde{\alpha}_i}{2} - 1, \frac{( L^{-1} \widetilde{U}
(L^T)^{-1} )_{ii}}{2} )$ and  $\operatorname{IG}( \frac{\widetilde{\alpha}_i}{2} - 1,
1 )$ are mutually absolutely continuous and the corresponding densities
with respect to Lebesgue measure are positive everywhere on $(0,
\infty)$ for all $i = 1,2, \ldots, m$. Also, since ${\bolds
\Theta}_{\psi, d_1, d_2}$ is bounded and closed, $\prod_{v=1}^{m-1}
\phi ( L^G_{\cdot v} \mid {\mathbf 0}, I_{n_v} ) \prod_{i=1}^m f_{\operatorname{IG}} (
\frac{\widetilde{\alpha}_i}{2} - 1, 1 )$ is bounded on ${\bolds
\Theta}_{\psi, d_1, d_2}$. Combining this with Proposition
\ref{adsgscndtn}, all required conditions in \cite{adscnvcmcs}, Theorem
6 are satisfied. Hence, the block Gibbs sampling Markov
chain, based on the derived conditional distributions, converges to the
desired stationary distribution $\pi_{U,\bolds{\alpha}}$.

 We note that in \cite{ghrmnslvgs}, page 18, the authors
introduce a procedure to sample from the $\mathcal{G}$-inverse Wishart
distributions (these are a narrow subclass of our priors $\pi^{P_G}_{U,
\bolds{\alpha}}$). Essentially, at every iteration, they cycle through
all of the rows of $\Sigma$. At the $i$th step in an iteration, they
sample the vector $\Sigma^G_{i \cdot} := (\Sigma_{ij})_{j \in
\mathcal{N} (i)}$ from its conditional distribution (Gaussian) given
the other entries of $\Sigma$ and then sample $\gamma_i :=
\frac{1}{\Sigma^{-1}_{ii}}$ from its conditional distribution
(inverse-gamma) given the other entries of $\Sigma$. Since $\Sigma$ is
a symmetric matrix, for $(i,j) \in E$, the variable $\Sigma_{ij}$
appears in $\Sigma^G_{i \cdot}$ as well as $\Sigma^G_{j \cdot}$. Hence,
$( ( \Sigma^G_{1 \cdot}, \gamma_1 ), ( \Sigma^G_{2 \cdot}, \gamma_2 ),
\ldots, ( \Sigma^G_{m \cdot}, \gamma_m ) )$ is not a disjoint partition
of the variable space. Therefore, their procedure is not strictly a
Gibbs sampling procedure and its convergence properties are not clear.
On the other hand, in our procedure, we cycle through $( L^G_{\cdot 1},
L^G_{\cdot 2}, \ldots, L^G_{\cdot m}, D )$, which is a disjoint
partition of the variable space. Hence, our procedure is a Gibbs
sampler in the true sense. There are also other differences between the
two procedures, such as the fact that $\gamma_i \neq D_{ii}$ unless $i=
m$.
\begin{remark*}
It is useful to compare our covariance priors to the conditionally
conjugate priors introduced by \cite{danlspoumi} in the complete case.
Upon closer investigation we discover that the priors of
\cite{danlspoumi} are quite different from $\widetilde{\pi}_{U, \bolds{
\alpha}} (L,D)$. First, they do not consider structural zeros. More
importantly, under their posterior, the distribution of $L^{-1}$
conditional on $D$ is jointly multivariate normal. In the general
decomposable covariance graph setting however, the zeros of $L$ do not
carry over to $L^{-1}$, and so it is not possible in our framework for
the distribution of the (constrained or unconstrained) elements of
$L^{-1}$, conditional on $D$, to be jointly multivariate normal.
\end{remark*}

\section{The special case of homogeneous graphs: Closed form expressions} \label{hprmrkvpty}

Note that the covariance graph model, that is, the family of
distributions
\[
\mathcal{G} = \{\mathcal{N}_m ({\mathbf 0}, \Sigma)\dvtx  \Sigma \in P_G\}
\]
(supported on $\mathbb{R}^m$) is a curved exponential family for any
connected noncomplete graph $G$. As discussed earlier, the fact that
the family is curved renders the Diaconis--Ylvisaker framework no
longer applicable in this setting. Hence, a rich and flexible class of
distributions was introduced in order to serve as priors for the class
of covariance graph models. A natural question to ask is whether the
class of priors itself belongs to a curved exponential
family.\footnote{Note: not the class of distributions associated with
the covariance graph probability model but rather the class of priors
that is introduced in this paper.} Indeed, this class of priors is
interesting in its own right and warrants an independent investigation.
Such analysis has the potential to place the class of priors in a known
framework and thus exploit this property.

 Let us therefore now turn our attention to the class of priors
$\{\widetilde{\pi}^{P_G}_{U,\bolds{\alpha}}\}_{U \in M_m^+}$ as a
family of distributions supported on $P_G$, with $U$ as a parameter. We
now state a lemma which formally establishes that the class of priors
can be framed in the context of natural exponential families.

\begin{lemma}
For arbitrarily fixed $\bolds{\alpha}$, the family of distributions
$\{\widetilde{\pi}^{P_G}_{U,\bolds{\alpha}}\}_{U \in M_m^+}$ is a
general exponential family, that is, it can be transformed into a
natural exponential family. The natural parameter is $U = ((U_{ij}))_{1
\leq i \leq j \leq m}$, the corresponding set of sufficient statistics
is $\Sigma^{-1} = ((\Sigma^{-1}_{ij}))_{1 \leq i \leq j \leq m}$ and
the cumulant generating function is $\log z_G (U,\bolds{\alpha})$.
\end{lemma}

\begin{pf}
The proof is straightforward and is therefore omitted.
\end{pf}

Placing the class of covariance priors in a natural exponential family
framework yields insights into the structure and functional form of
this class of distributions. As noted earlier, $z_G (U,
\bolds{\alpha})$ is not generally available in closed form. A question that
naturally arises is whether there are any conditions under which $z_G
(U,\bolds{\alpha})$ can be evaluated in closed form. In this section,
we establish that when $G$ is homogeneous, $z_G (U, \bolds{\alpha})$
and ${\mathbf E}_{U,  \bolds{\alpha}}$ can be evaluated in closed form.
It is known that when $G$ is a homogeneous graph, the covariance graph
model is Markov equivalent to an appropriate DAG (see
\cite{drtnrchrn3}). It is, however, important to clarify that the
Markov equivalence of covariance graph models and DAGs does not
immediately imply that Bayesian inference for covariance graph models
using our priors automatically follows. We also need to establish a
correspondence between our priors and known priors for DAG models. We
now prove that in the special case when $G$ is homogeneous, our priors
correspond to the standard conjugate priors for an appropriate DAG.
This yields yet another property of our class of priors. The following
theorem is the main result of this section and helps us establish the
aforementioned correspondence.

\begin{theorem} \label{trnfmdcpic}
Let $G = (V,E)$ be homogeneous, with vertices ordered according to the
Hasse perfect vertex elimination scheme specified in Section
\ref{hmgnsgphpy}, that is, $V \in S_H$. If $\Sigma \sim \pi^{P_G}_{U,
\alpha}$ and $\Sigma = LDL^T$ is its modified Cholesky
decomposition, then
\[
\{ ( D_{ii}, (\Sigma^{\prec i})^{-1} \Sigma^\prec_{\cdot i} ) \}_{1
\leq i \leq m}
\]
are mutually independent. Furthermore, the distributions of these
quantities are specified as follows:
%
\begin{eqnarray}
(\Sigma^{\prec i})^{-1} \Sigma^\prec_{\cdot i} \mid D_{ii} &\sim&
\mathcal{N} ((U^{\prec i})^{-1} U^\prec_{\cdot i}, D_{ii} (U^{\prec i})^{-1} );\nonumber\\
D_{ii} &\sim& \operatorname{IG} \biggl( \frac{\alpha_i}{2} - \frac{|\mathcal{N}^\prec
(i)|}{2} - 1, \frac{U_{ii} - (U^\prec_{\cdot i})^T (U^{\prec i})^{-1}
U^\prec_{\cdot i}}{2}\biggr)\nonumber\\
\eqntext{\forall i = 1,2, \ldots, m.}
\end{eqnarray}
\end{theorem}

\begin{remark*}
The above result decomposes $\Sigma$ into mutually
independent coordinates. Note that for any $i$ such that $\bar{i}$ is a
leaf of the Hasse tree and $i$ has the minimal label in its equivalence
class $\bar{i}$, we have
\[
\mathcal{N}^\prec (i) = \phi.
\]
In this case, it is understood that $\Sigma^{\prec i}$ and
$\Sigma^\prec_{\cdot i}$ are vacuous parameters and that $D_{ii} =
\Sigma_{ii}$.
\end{remark*}

\begin{pf*}{Proof of Theorem \ref{trnfmdcpic}}
Let $G$ be a homogeneous graph with $m$ vertices, with the
vertices ordered according to the Hasse perfect elimination scheme
specified in Section \ref{hmgnsgphpy}. Recall that the vertices of the
Hasse diagram of $G$ are equivalence classes formed by the relation $R$
defined in Section \ref{hmgnsgphpy}. The vertex labeled $m$ clearly
lies in the equivalence class of vertices at the root of the
corresponding Hasse diagram. Let us remove the vertex labeled $m$ from
the graph $G$ and let $G'$ denote the induced graph on the remaining $m
- 1$ vertices. The graph $G'$ can be of the following two types.
\begin{itemize}
\item Case I: If the equivalence class of $m$ contains more than one
element, then $G'$ is a homogeneous graph with the Hasse diagram having
the same depth as the\vadjust{\goodbreak} Hasse diagram of $G$, but with one less vertex in
the equivalence class at the root. Recall that the depth of a tree is
the length of the longest path from its root to any leaf.
\item Case II: If the equivalence class of $m$ contains only one element, then
$G'$ is a disconnected graph, with the connected components being
homogeneous graphs with the Hasse diagram having depth one less than
the depth of the Hasse diagram of $G$.
\end{itemize}

Note that for every $1 \leq i \leq m$ such that $\mathcal{N}^\prec (i)
\neq \phi,  \Sigma^{\preceq i}$ can be partitioned as
\[
\Sigma^{\preceq i} = \left[\matrix{\Sigma^{\prec i} &
\Sigma^\prec_{\cdot i}\vspace*{2pt} \cr (\Sigma^\prec_{\cdot i})^T & \Sigma_{ii}}\right].
\]
Also, note that if $Z \sim \mathcal{N} ({\mathbf 0}, \Sigma),$ then
$D_{ii}$ is the conditional variance of $Z_i$ given $Z_1, Z_2, \ldots,
Z_{i-1}$ (see \cite{hngluprmlu}). Note that $\Sigma_{kl} = 0$ for all
$1 \leq k,l \leq i,  k \in \mathcal{N}^\prec (i),  l \notin
\mathcal{N}^\prec (i)$. It follows that $D_{ii} = \Sigma_{ii} -
(\Sigma^\prec_{\cdot i})^T (\Sigma^{\prec i})^{-1} \Sigma^\prec_{\cdot
i}$. Hence, by the formula for the inverse of a partitioned matrix, it
follows that
\[
(\Sigma^{\preceq i})^{-1} = \left[\matrix{\displaystyle(\Sigma^{\prec i})^{-1} +
\frac{( (\Sigma^{\prec i})^{-1} \Sigma^\prec_{\cdot i} ) (
(\Sigma^{\prec i})^{-1} \Sigma^\prec_{\cdot i} )^T}{D_{ii}} &
\displaystyle-\frac{(\Sigma^{\prec i})^{-1} \Sigma^\prec_{\cdot
i}}{D_{ii}}\vspace*{2pt}
\cr \displaystyle-\frac{( (\Sigma^{\prec i})^{-1} \Sigma^\prec_{\cdot
i} )^T}{D_{ii}} & \dfrac{1}{D_{ii}}}\right].
\]
Hence,
\begin{eqnarray}\label{mrkvdcmptr}
&&\operatorname{tr} ( (\Sigma^{\preceq i})^{-1} U^{\preceq i} )\nonumber\\
&&\qquad= \operatorname{tr} ( (\Sigma^{\prec i})^{-1} U^{\prec i} )
\nonumber\\ [-8pt]\\ [-8pt]
&&\qquad\quad{}+ \frac{1}{D_{ii}} \bigl( (\Sigma^{\prec i})^{-1} \Sigma^\prec_{\cdot i} -
(U^{\prec i})^{-1} U^\prec_{\cdot i} \bigr)^T U^{\prec i} \bigl( (\Sigma^{\prec i})^{-1}
\Sigma^\prec_{\cdot i} - (U^{\prec i})^{-1} U^\prec_{\cdot i} \bigr)  \nonumber\\
&&\qquad\quad{} +\frac{1}{D_{ii}} \bigl( U_{ii} - (U^\prec_{\cdot i})^T (U^{\prec i})^{-1} U^\prec_{\cdot i}
\bigr).\nonumber
\end{eqnarray}
We again note that from our argument at the beginning of the proof,
$\Sigma^{\prec i} = \Sigma^{\preceq (i-1)}$ or $\Sigma^{\prec i}$ has a
block diagonal structure (after an appropriate permutation of the rows
and columns) with blocks $\Sigma^{\preceq i_1}, \Sigma^{\preceq i_2},
\ldots, \Sigma^{\preceq i_k}$  for some $k > 1,  1 \leq i_1, i_2,
\ldots, i_k < i$.\footnote{If the equivalence class of $i$ has $k$
children in the Hasse diagram of $G$ and $V_j$  is the set of vertices
in $V$ belonging to the subtree rooted at the $j$th child, then $V_j$,
for $1 \leq j \leq k$, are disjoint subsets. In fact, if $i_j = \max
\{i'\dvtx  i' \in V_j\}$, then it follows by the construction of the Hasse
diagram that $V_j = \mathcal{N}^\preceq (i_j)$ for $1 \leq j \leq k$.}
It follows that
\[
\operatorname{tr} ( (\Sigma^{\prec i})^{-1} U^{\prec i} ) = \sum_{j=1}^k \operatorname{tr} (
(\Sigma^{\preceq i_j})^{-1} U^{\preceq i_j} ).
\]
Note that $\Sigma = \Sigma^{\preceq m}$. Using (\ref{mrkvdcmptr})
recursively, we get
\begin{eqnarray}\label{traceidnty}
& & \operatorname{tr}(\Sigma^{-1} U) \nonumber\\
&&\qquad= \sum_{i=1}^m \frac{1}{D_{ii}} \bigl( (\Sigma^{\prec i})^{-1}
\Sigma^\prec_{\cdot i} - (U^{\prec i})^{-1} U^\prec_{\cdot i} \bigr)^T
U^{\prec i} \bigl( (\Sigma^{\prec i})^{-1}
\Sigma^\prec_{\cdot i} - (U^{\prec i})^{-1} U^\prec_{\cdot i} \bigr) \\
&&\qquad\quad\hphantom{\sum_{i=1}^m}{} + \frac{1}{D_{ii}} \bigl( U_{ii} - (U^\prec_{\cdot i})^T (U^{\prec
i})^{-1} U^\prec_{\cdot i} \bigr).\nonumber
\end{eqnarray}
Let us now evaluate the Jacobian of the transformation
\[
\Sigma \rightarrow \{ ( D_{ii}, (\Sigma^{\prec i})^{-1}
\Sigma^\prec_{\cdot i} ) \}_{1 \leq i \leq m}.
\]
It follows by simple matrix manipulations that  the Jacobian of the
transformation
\[
\Sigma^{\preceq i} \rightarrow ( \Sigma^{\prec i}, (\Sigma^{\prec
i})^{-1} \Sigma^\prec_{\cdot i}, D_{ii} )
\]
is given by $| \Sigma^{\prec i} |$.

 Once more, note that $\Sigma = \Sigma^{\preceq m}$ and, as mentioned earlier, $\Sigma^{\prec i} =
\Sigma^{\preceq (i-1)}$ or $\Sigma^{\prec i}$ (after an appropriate
permutation of the rows and columns) has a block diagonal structure
with blocks $\Sigma^{\preceq i_1}, \Sigma^{\preceq i_2}, \ldots,
\Sigma^{\preceq i_k}$ for some $k > 1,  1 \leq i_1, i_2, \ldots, i_k <
i$. Hence, by regarding the transformation
\[
\Sigma \rightarrow \{ ( D_{ii}, (\Sigma^{\prec i})^{-1}
\Sigma^\prec_{\cdot i} ) \}_{1 \leq i \leq m}
\]
as a series of transformations of the type $\Sigma^{\preceq i}
\rightarrow ( \Sigma^{\prec i}, (\Sigma^{\prec i})^{-1}
\Sigma^\prec_{\cdot i}, D_{ii} )$, it follows that the determinant of
the Jacobian is given by
%
\begin{equation} \label{jcbnhgsgph}
\prod_{i=1}^m | \Sigma^{\prec i} | = \prod_{i=1}^m \prod_{j \in
\mathcal{N}^\prec(i)} D_{jj} =  \prod_{j=1}^m D_{jj}^{n_j}.
\end{equation}
Here, as in Section \ref{inducedprr} Lemma \ref{jcbnttpqtp},
\[
n_j = |\{(i>j\dvtx (i,j) \in E\}|\qquad\forall j = 1,2, \ldots, m.
\]
Also, from Section \ref{inducedprr},
\[
\pi^{P_G}_{U,\bolds{\alpha}} (\Sigma) = \frac{1}{z_G (U,
\bolds{\alpha})} e^{-( \operatorname{tr}( \Sigma^{-1} U ) + \sum_{j=1}^m (2 n_j +
\alpha_j) \log D_{jj} )/2},\qquad\Sigma \in P_G.
\]
Let
\[
{\bolds \Gamma} = \{ ( D_{ii}, (\Sigma^{\prec i})^{-1}
\Sigma^\prec_{\cdot i} ) \}_{1 \leq i \leq m}.
\]
It follows from the decomposition of $\operatorname{tr}(\Sigma^{-1} U )$ from
(\ref{traceidnty}) and the computation of the determinant of the
Jacobian (\ref{jcbnhgsgph}) that
%
\begin{eqnarray}\label{hmgnsdtyfr}
&& \pi^{\Gamma}_{U, \bolds{\alpha}} ( \{ ( D_{ii}, (\Sigma^{\prec
i})^{-1}\Sigma^\prec_{\cdot i} )_{1 \leq i \leq m} \} ) \nonumber\\
&&\qquad= \frac{1}{z_G (U,\bolds{\alpha})} \prod_{i=1}^m e^{-
1/(2D_{ii}) ( (\Sigma^{\prec i})^{-1} \Sigma^\prec_{\cdot i} -
(U^{\prec i})^{-1} U^\prec_{\cdot i} )^T U^{\prec i} ( (\Sigma^{\prec i})^{-1}
\Sigma^\prec_{\cdot i} - (U^{\prec i})^{-1} U^\prec_{\cdot i} )}\hspace*{-16pt} \\
&&\qquad\quad{}\times \prod_{i=1}^m e^{- 1/(2D_{ii}) ( U_{ii}
- (U^\prec_{\cdot i})^T (U^{\prec i})^{-1} U^\prec_{\cdot i} )}
D_{ii}^{-\alpha_i/2}.\nonumber
\end{eqnarray}
The above proves the mutual independence of $\{ ( D_{ii},
(\Sigma^{\prec i})^{-1} \Sigma^\prec_{\cdot i} ) \}_{1 \leq i \leq m}$.
From the joint density of $( D_{ii}, (\Sigma^{\prec i})^{-1}
\Sigma^\prec_{\cdot i} )$, it is clear that
\[
(\Sigma^{\prec i})^{-1} \Sigma^\prec_{\cdot i} \mid D_{ii} \sim
\mathcal{N} ( (U^{\prec i})^{-1} U^\prec_{\cdot i}, D_{ii} (U^{\prec
i})^{-1} ).
\]
To evaluate the marginal density of $D_{ii}$, we integrate out
$(\Sigma^{\prec i})^{-1} \Sigma^\prec_{\cdot i}$ from the joint density
of $(D_{ii}, (\Sigma^{\prec i})^{-1} \Sigma^\prec_{\cdot i})$. Note
that
{\fontsize{10.84}{12}\selectfont{
\begin{eqnarray*}
\hspace*{-4pt}&&\int_{\mathbb{R}^{|\mathcal{N}^\prec (i)|}}e^{- 1/(2D_{ii}) (
(\Sigma^{\prec i})^{-1} \Sigma^\prec_{\cdot i} - (U^{\prec i})^{-1}
U^\prec_{\cdot i} )^T U^{\prec i} ( (\Sigma^{\prec i})^{-1}
\Sigma^\prec_{\cdot i} - (U^{\prec i})^{-1} U^\prec_{\cdot i} )}\, d
((\Sigma^{\prec i})^{-1} \Sigma^\prec_{\cdot i})\\
\hspace*{-4pt}&&\qquad = C
D_{ii}^{-|\mathcal{N}^\prec (i)|/2},
\end{eqnarray*}}}%
where $C$ is a constant, since the above integral is essentially an
unnormalized multivariate normal integral. Hence, the marginal density
of $D_{ii}$ is given by
\[
\pi^{D_{ii}}_{U,\bolds{\alpha}} (d) \propto e^{-( U_{ii} -
(U^\prec_{\cdot i})^T (U^{\prec i})^{-1} U^\prec_{\cdot i} )/(2d)}
d^{-( \alpha_{i}/2 + |\mathcal{N}^\prec (i)|/2 )}.
\]
We can therefore conclude that
\[
D_{ii} \sim \operatorname{IG} \biggl( \frac{\alpha_i}{2} - \frac{|\mathcal{N}^\prec (i)|}{2}
- 1, \frac{U_{ii} - (U^\prec_{\cdot i})^T (U^{\prec i})^{-1}
U^\prec_{\cdot i}}{2} \biggr).
\]
\upqed\end{pf*}

\begin{remark*}
At  first glance, it seems as if the only part of $U$
that appears in Theorem \ref{trnfmdcpic} is $(U_{ij})_{(i,j) \in E}$,
that is, the projection of $U$ onto $I_G$. Hence, one could incorrectly
conclude that up to the number of shape parameters in each equivalence
class, in the homogeneous case, the priors introduced in this paper and
the $\mathit{IW}_{Q_G}$ are identical. However, a careful inspection shows that
this is not the case. Note that the conditional covariance of
$(\Sigma^{\prec i})^{-1} \Sigma^{\prec}_{\cdot i}$ is $D_{ii} (U^{\prec
i})^{-1}$, and $U^{\prec i}$ can contain entries of the form $U_{kl}$
such that $(k,l) \notin E$.  For example, suppose that $G =
\stackrel{1}{\bullet} - \stackrel{3}{\bullet} -
\stackrel{2}{\bullet} $. Then,
\[
 U^{\prec 3} = \pmatrix{
 U_{11} & U_{12} \cr
 U_{21} & U_{22}
},
\]
but $(1,2) \notin E$. Hence, in the homogeneous setting, $\pi^{P_G}_{U,
\bolds{\alpha}}$ is truly a larger class than the $\mathit{IW}_{Q_G}$ family of
distributions.

We now establish the correspondence between $\pi^\Gamma_{U,
\bolds{\alpha}}$ in (\ref{hmgnsdtyfr}) and the conjugate prior for an
appropriate Gaussian DAG model. Let $G = (V,E)$ be a homogeneous graph
with $V \in S_H$, that is, the vertices have been ordered according to
the perfect vertex elimination scheme for homogeneous graphs outlined
in Section \ref{hmgnsgphpy}. Let us construct a DAG as follows:
\begin{enumerate}
\item Consider the Hasse diagram of $G$ (for simplicity and clarity of
exposition, assume that the equivalence class at each vertex has just
one element).
\item Assign a directed edge from $u$ to $v$ if $u$ is a
descendant of $v$ in the Hasse tree, that is, reverse the directions of
all the arrows, including those that do not appear in the Hasse tree,
but which are implied by transitivity.
\end{enumerate}

An example of a DAG constructed in this manner is given in Figure
\ref{exmlhgsfig}.
%
\begin{figure}

\includegraphics{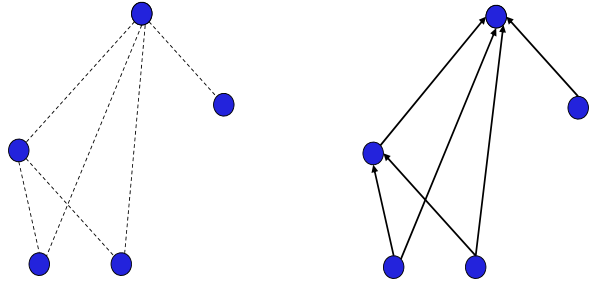}

\caption{Example of homogeneous covariance graph (left);
associated DAG model (right).} \label{exmlhgsfig}
\end{figure}

Now, let $\operatorname{pa}(i)$ denote the set of parents of $i$ according to the
direction specified above. If ${\mathbf Y} \sim \mathcal{N}_m ({\mathbf
0}, \Sigma)$ has a distribution which is Markov with respect to the
above DAG, then the density of ${\mathbf Y}$ factorizes as
\begin{eqnarray*}
f({\mathbf y}) &=& \prod_{i=1}^m f\bigl(y_i \mid {\mathbf y}_{\operatorname{pa}(i)}\bigr)\\
&=& \prod_{i=1}^m \frac{1}{\sqrt{2 \pi} D_{ii}^{1/2}} e^{- (
y_i - ( (\Sigma^{\prec i})^{-1} \Sigma^\prec_{\cdot i} )^T {\mathbf
y}_{\operatorname{pa}(i)} )^T D_{ii}^{-1} ( y_i - ( (\Sigma^{\prec i})^{-1}
\Sigma^\prec_{\cdot i} )^T {\mathbf y}_{\operatorname{pa}(i)} )},
\end{eqnarray*}
where $D_{ii} := \Sigma_{ii} - (\Sigma^\prec_{\cdot i})^T
(\Sigma^{\prec i})^{-1} \Sigma^\prec_{\cdot i} = \Sigma_{i \mid \operatorname{pa} (i)}
$.

The standard conjugate prior for each factor of the product above can
be obtained as follows. Given an arbitrary positive definite matrix $U$
and $\alpha'_i > 0$ for $i = 1,2, \ldots, m$, let
\begin{eqnarray*}
(\Sigma^{\prec i})^{-1} \Sigma^\prec_{\cdot i} \mid D_{ii} &\sim&
\mathcal{N} ( (U^{\prec i})^{-1} U^\prec_{\cdot i},
D_{ii}(U^{\prec i})^{-1} ),\\
 D_{ii} &\sim& \operatorname{IG} \biggl( \alpha'_i, \frac{U_{ii} - (U^\prec_{\cdot i})^T
(U^{\prec i})^{-1} U^\prec_{\cdot i}}{2}\biggr)
\end{eqnarray*}
for $i = 1,2, \ldots, m$, where $\{ ( D_{ii}, (\Sigma^{\prec i})^{-1}
\Sigma^\prec_{\cdot i} ) \}_{1 \leq i \leq m}$ are mutually
independent. This corresponds precisely to the $\pi^\Gamma_{U, {\bolds{\alpha}}}$ density in (\ref{hmgnsdtyfr}).

We now proceed to state, without proof, results for homogeneous
covariance graph models by exploiting the correspondence of our priors
to the standard conjugate priors for DAGs. In particular, hyper-Markov
properties, the normalizing constant and expected values for covariance
graph models, for $G$ homogeneous, are formally stated below.

Let $G = (V,E)$ be a homogeneous graph with $V \in S_H$, that is, the
vertices have been ordered according to the perfect vertex elimination
scheme for homogeneous graphs outlined in Section \ref{hmgnsgphpy}. Let
$D$ be the directed graph obtained from $G$ by directing all edges in
$G$ from the vertex with the smallest number to the vertex with the
highest number. Let $\operatorname{pa} (i)$ denote the set of parents of $i$ according
to the direction specified in $D$. It follows that $\operatorname{pa} (i) =
\mathcal{N}^\prec (i)$. As in \cite{dawidlrtzn,ltcmssmwdg}, let $\operatorname{pr} (i)
= \{1,2, \ldots, i-1\}$ denote the set of predecessors of $i$ according
to the direction specified in $D$. We now proceed to define the
hyper-Markov property.
\end{remark*}

\begin{definition}
A family of priors $\mathcal{F}$ on $P_G$ satisfies the strong
hyper-Markov property with respect to the direction $D$ if, whenever
$\pi \in \mathcal{F}$ and $\Sigma \sim \pi$,
\[
\Sigma_{i \mid \operatorname{pa} (i)} \perp \Sigma_{\operatorname{pr} (i)}\qquad\forall  1 \leq i \leq m,
\]
where $\Sigma_{i \mid \operatorname{pa} (i)} := \Sigma_{ii} - (\Sigma^\prec_{\cdot
i})^T (\Sigma^{\prec i})^{-1} \Sigma^\prec_{\cdot i} = D_{ii}$.
\end{definition}

In the following corollary, we state, without proof,  that the family
of priors $\pi^{P_G}_{U, \bolds{\alpha}}$ satisfies the strong
hyper-Markov property with respect to the direction $D$.

\begin{cor}\label{hprmrkvstm}
Let $G = (V,E)$ be homogeneous with $V \in S_H$. If $\Sigma \sim
\pi^{P_G}_{U,\bolds{\alpha}}$, then
\[
D_{ii} \perp \Sigma_{\{1,2, \ldots, i-1\}}\qquad\forall  1 \leq i \leq m.
\]
\end{cor}

\begin{remark*} Recall that
\[
\Sigma^{\prec i} = ((\Sigma_{uv}))_{u,v \in \mathcal{N}^\prec (i)}
\]
is different from
\[
\Sigma_{\{1,2, \ldots, i-1\}} = ((\Sigma_{uv}))_{1 \leq u,v \leq i-1}.
\]
\end{remark*}

We demonstrated in Section \ref{ltcmssmcpn} that the family $\mathit{IW}_{Q_G}$
of Letac and Massam \cite{ltcmssmwdg} is a subfamily of our class of
priors $\pi^{P_G}_{U,\bolds{\alpha}}$ when $G$ is homogeneous.
Consequently, we can now prove hyper-Markov properties for the
$\mathit{IW}_{Q_G}$ family.

\begin{cor}
Let $G = (V,E)$ be homogeneous with $V \in S_H$. Let $D$ be the
directed graph obtained from $G$ by directing all edges in $G$ from the
vertex with the smallest number to the vertex with the highest number.
The family $\mathit{IW}_{Q_G}$ is then strong hyper-Markov with respect to the
direction $G_H$.
\end{cor}

Hyper-Markov properties for the $\mathit{IW}_{Q_G}$ family were not established
in \cite{ltcmssmwdg}. Hence, we note that the corollary above is a new
result for this family.

We now proceed to state, without proof, the functional form of the
normalizing constant for homogeneous graphs, once again exploiting the
correspondence between our covariance priors and the conjugate priors
for DAGs. In particular, below, we state  necessary and sufficient
conditions for existence of the normalizing constant and give an
explicit expression for it in such cases.

\begin{cor} \label{nrmlcnsthn}
Let $G = (V,E)$ be a homogeneous graph with vertices ordered such that
$V \in S_H$. Then, $z_G (U,\bolds{\alpha}) < \infty$ if and only if
$\bolds{\alpha}$ satisfies the conditions in Proposition
\ref{nrmcnstdml}, that is, $\alpha_i > |\mathcal{N}^\prec (i)| + 2\
\forall i = 1,2, \ldots, m$. In this case,
%
\begin{eqnarray} \label{nrmlcnsten}
&&z_G (U,\bolds{\alpha}) = \prod_{i=1}^m\Gamma \biggl(
\frac{\alpha_i}{2} - \frac{|\mathcal{N}^\prec (i)|}{2} - 1 \biggr)
2^{\alpha_i/2 - 1} (\sqrt{\pi})^{|\mathcal{N}^\prec (i)|}\nonumber \\
[-8pt]\\ [-8pt]
&&\hphantom{z_G (U,\alpha) = \prod_{i=1}^m}{}\times|U^{\prec i} |^{\alpha_i/2 - |\mathcal{N}^\prec (i)|/2 -
3/2}\big/| U^{\preceq i} |^{\alpha_i/2 -
|\mathcal{N}^\prec (i)|/2 - 1}.\nonumber
\end{eqnarray}
\end{cor}

We now proceed to state, without proof, expected values related to our
class of priors $\pi^{P_G}_{U,\bolds{\alpha}}$ when $G$ is
homogeneous, again by exploiting the correspondence between our
covariance priors and conjugate priors for DAGs. In particular, we now
provide a recursive method that gives closed form expressions for the
expected value of the covariance matrix when $\Sigma \sim \pi^{P_G}_{U,
\bolds{\alpha}}$. Since $\Sigma_{uv} = 0\ \forall (u,v) \notin E$, we
only need to evaluate the expectation of $\Sigma_{ii}$ and
$\Sigma^\prec_{\cdot i}$ for every $1 \leq i \leq m$. Let
\[
A_1 := \{i \in V\dvtx \mathcal{N}^\prec (i) = \phi\}.
\]

Clearly, if $i \in A_1$, then $\Sigma^{\prec i}$ and $\Sigma_{\cdot i}$
are vacuous parameters and $D_{ii} = \Sigma_{ii}$. It follows from
Theorem \ref{trnfmdcpic} that for $i \in A_1$,
\[
{\mathbf E}_{U,\bolds{\alpha}} [\Sigma_{ii}] = {\mathbf E}_{U,
\bolds{\alpha}} [D_{ii}] = \frac{U_{ii} - (U^\prec_{\cdot i})^T
(U^{\prec i})^{-1} U^\prec_{\cdot i}}{\alpha_i - 4},
\]
assuming that $\alpha_i > 4$, since $X \sim \operatorname{IG}(\lambda, \gamma)$
implies that ${\mathbf E}[X] = \frac{\lambda}{\gamma - 1}$.

For $k = 2,3,4, \ldots,$ define
\[
A_k = \Biggl\{ i \in V\dvtx   \mathcal{N}^\prec (i) \subseteq \bigcup_{l=1}^{k-1}
A_l \Biggr\} \bigg\backslash \Biggl( \bigcup_{l=1}^{k-1} A_l \Biggr).
\]
Since there are finitely many vertices in $V$, there exists some $k^*$
such that $A_k = \phi$ for $k > k^*$. The sets $\{A_k\}_{1 \leq k \leq
k^*}$ essentially provide a way of computing ${\mathbf E}_{U,
\bolds{\alpha}} [ \Sigma ]$, by starting at the bottom of the Hasse diagram of
$G$ and then moving up sequentially.

\begin{cor} \label{exptncovar}
Let G be a homogeneous graph. Given the expectations of
$\Sigma^\prec_{\cdot j}$ and $\Sigma_{jj}$ for $j \in \bigcup_{l=1}^{k-1}
A_l$, the expectations of $\Sigma^\prec_{\cdot i}$ and $\Sigma_{ii}$
for $i \in A_k$ are given, respectively, by the  expressions
\begin{eqnarray*}
{\mathbf E}_{U, \bolds{\alpha}} [\Sigma^\prec_{\cdot i}] &=&
{\mathbf E}_{U, \bolds{\alpha}} [\Sigma^{\prec i} ] (U^{\prec i})^{-1} U^\prec_{\cdot i};\\
{\mathbf E}_{U, \bolds{\alpha}} [\Sigma_{ii}] &=& \frac{U_{ii} -
(U^\prec_{\cdot i})^T(U^{\prec i})^{-1} U^\prec_{\cdot i}}{\alpha_i - |\mathcal{N}^\prec (i)| - 4} \\
& & {}+\operatorname{tr}\biggl( {\mathbf E}_{U,\bolds{\alpha}} [ \Sigma^{\prec i} ]
 \biggl( \frac{(U^{\prec i})^{-1} ( U_{ii} - (U^\prec_{\cdot i})^T (U^{\prec i})^{-1}
U^\prec_{\cdot i} )}{\alpha_i - |\mathcal{N}^\prec (i)| - 4}\\
&& \hspace*{108pt}{}+(U^{\prec i})^{-1} U^\prec_{\cdot i} (U^\prec_{\cdot i})^T (U^{\prec
i})^{-1} \biggr)\biggr ),
\end{eqnarray*}
provided that $\alpha_i > |\mathcal{N}^\prec (i)| + 4$.
\end{cor}

The corollary is not formally proved since it follows directly from the
correspondence between our covariance priors and the natural conjugate
priors for DAGs. We note once more that the expressions above yield a
recursive but closed form method to calculate ${\mathbf E} [\Sigma]$
when $\Sigma \sim \pi^{P_G}_{U,\bolds{\alpha}}$.

\begin{remark*}
There is an intriguing parallel between the expressions
for the normalizing constant and the expected values for the
$\pi^{P_G}_{U,\bolds{\alpha}}$ distribution and the $\mathit{IW}_{P_G}$
distribution (as derived in \cite{ltcmssmwdg,rajmasscar}) when $G$ is
homogeneous. This automatically leads one to wonder if the
$\pi^{P_G}_{U,\bolds{\alpha}}$ and the $\mathit{IW}_{P_G}$ distributions are
the same. We now show that this is not the case.

If one compares the density of $( D_{ii}, (\Sigma^{\prec i})^{-1}
\Sigma^\prec_{\cdot i} )_{1 \leq i \leq m}$ in (\ref{hmgnsdtyfr}) and the
$\mathit{IW}_{P_G}^*$ density in $(3.16)$ of \cite{ltcmssmwdg},
 they initially appear to have the same functional form.
 We now proceed to show that they are supported on different spaces.
 This difference is illustrated by the following example. Let
 $G =\,\stackrel{1}{\bullet} -
\stackrel{3}{\bullet} - \stackrel{2}{\bullet}$. In this case, each
equivalence class in the Hasse tree of $G$ has exactly one vertex. Note
that vertex $3$ has two descendants, and vertices $1$ and $2$ do not
have any descendants in the Hasse tree of $G$. Hence, it follows that
the density of $( D_{ii}, (\Sigma^{\prec i})^{-1} \Sigma^\prec_{\cdot
i} )_{1 \leq i \leq 3}$ is supported on $( \mathbb{R} \times \mathit{NULL},
\mathbb{R} \times \mathit{NULL}, \mathbb{R} \times \mathbb{R}^2 )$. On the other
hand, vertices $1$ and $2$ have one ancestor, and vertex $3$ has no
ancestors in the Hasse tree of $G$. Hence, it follows that the
$\mathit{IW}_{P_G}^*$ density is supported on $( \mathbb{R} \times \mathbb{R},
\mathbb{R} \times \mathbb{R}, \mathbb{R} \times \mathit{NULL} )$.

So, at first glance, it looks as if $( D_{ii}, (\Sigma^{\prec i})^{-1}
\Sigma^\prec_{\cdot i} )_{1 \leq i \leq m}$ has the same form as
$\mathit{IW}_{P_G}^*$, but, upon further examination, we see that even for the
simplest homogeneous graph, they are structurally different. In fact, $
(\mathbb{R} \times \mathit{NULL}, \mathbb{R} \times \mathit{NULL}, \mathbb{R} \times
\mathbb{R}^2 ) $ does not support the $\mathit{IW}_{P_G}^*$ distribution for any
$G$ that is homogeneous.
\end{remark*}

\section{Examples} \label{examplessection}

The main purpose of this paper is to undertake a theoretical
investigation of our class of distributions and their efficacy for use
in Bayesian estimation in covariance graph models. We nevertheless
provide two examples (one real and one simulated) to demonstrate how
the methodology developed in this paper can be implemented.

\subsection{Genomics example}

We provide an illustration of our methods on a data set consisting of
gene expression data from microarray experiments with yeast strands
from Gasch et al. \cite{gskcesbbyd}. This data set has also been
analyzed in \cite{chaudrtric,drtnrchrn3}. As in \cite{chaudrtric,drtnrchrn3},
we consider a subset of eight genes involved in galactose
utilization. There are $n = 134$ experiments and the empirical
covariance matrix for these measurements is provided in Table
\ref{ystcvncmtx}. Note that the sample covariance matrix is obtained
after centering since the mean is not assumed to be zero.

\begin{table}
\tabcolsep=0pt
\caption{Empirical covariance matrix for yeast data}\label{ystcvncmtx}
\begin{tabular*}{\textwidth}{@{\extracolsep{\fill}}lcccccccc@{}}
\hline
 & \textbf{GAL11} & \textbf{GAL4} & \textbf{GAL80} & \textbf{GAL3} & \textbf{GAL7} & \textbf{GAL10} & \textbf{GAL1} & \textbf{GAL2}\\
\hline
GAL11 & 0.152 & & & & & & & \\
GAL4 & 0.034 & 0.130 & & & & & & \\
GAL80 & 0.015 & 0.039 & 0.221 & & & & & \\
GAL3 & $-$0.055\phantom{$-$} & 0.034 & 0.073 & 0.608 & & & & \\
GAL7 & $-$0.051\phantom{$-$} & $-$0.053\phantom{$-$} & 0.183 & 0.722 & 3.423 & & & \\
GAL10 & $-$0.048\phantom{$-$} & $-$0.039\phantom{$-$} & $-$0.188\phantom{$-$} & 0.553 & 2.503 & 2.372 & & \\
GAL1 & $-$0.066\phantom{$-$} & $-$0.061\phantom{$-$} & 0.224 & 0.517 & 2.768 & 2.409 & 2.890 & \\
GAL2 & $-$0.119\phantom{$-$} & $-$0.018\phantom{$-$} & 0.208 & 0.583 & 2.547 & 2.278 & 2.514 &
2.890 \\
\hline
\end{tabular*}
\end{table}

\begin{figure}[b]

\includegraphics{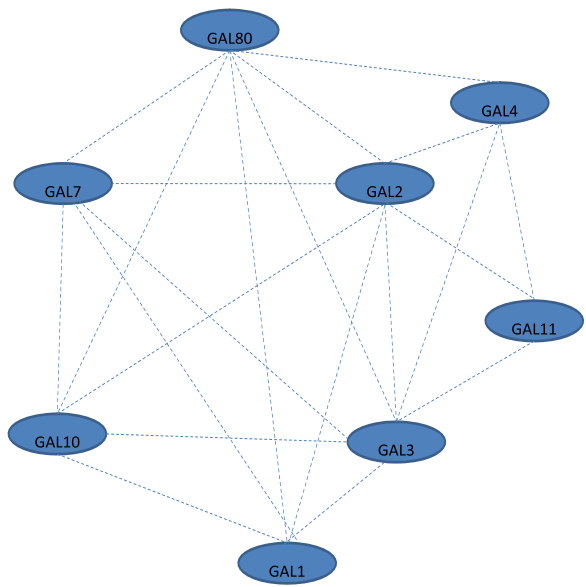}

\caption{Covariance graph for yeast data.} \label{ystexpngph}
\end{figure}

 We consider the covariance graph model specified by the graph
$G$ in Figure \ref{ystexpngph} with the overall aim of estimating
$\Sigma$ under this covariance graph model. The maximum likelihood
estimate for $\Sigma \in P_G$, provided by the iterative conditional
fitting algorithm described in \cite{chaudrtric}, yields a deviance of
$4.694$ over $7$ degrees of freedom, thus indicating a good model fit.
The
maximum likelihood estimate is provided in Table \ref{ystcvncmtx}. We use the following ordering for our analysis:
$\{\mathrm{GAL}11, \mathrm{GAL}4, \mathrm{GAL}80, \mathrm{GAL}3, \mathrm{GAL}7, \mathrm{GAL}10, \mathrm{GAL}1, \mathrm{GAL}2\}$.

Our goal is to obtain the posterior mean for $\Sigma$ under our
new class of priors and then to provide Bayes estimators for $\Sigma$.
We use two diffuse priors to illustrate our methodology. The first
prior is denoted as $\widetilde{\pi}_{U_1, \bolds{\alpha}^1}$, where
$U_1 = \frac{\operatorname{tr}(S)}{8} I_8,  \alpha^1_i = 5 + |\mathcal{N}^\prec (i)|,
 i = 1,2, \ldots, 8$, that is, $ \bolds{\alpha}_1 =
(5,6,6,8,7,8,9,12)$. The second prior used is $\widetilde{\pi}_{U_2,
\bolds{\alpha}^2}$, where $U_2 = 0$, $\alpha^2_i = 2$, $i = 1,2,
\ldots, 8$. Note that we could have used any ordering in $S_D$ for our
analysis. As an example, we select an alternate ordering, $\{\mathrm{GAL}11,
\mathrm{GAL}4, \mathrm{GAL}80, \mathrm{GAL}10, \mathrm{GAL}2, \mathrm{GAL}3, \mathrm{GAL}1,\break \mathrm{GAL}7\}$, and also consider the
two priors mentioned above under this alternative ordering. The block
Gibbs sampling procedure was run for the four priors as specified in
Section~\ref{cndtlcmptn}. The burn-in period was chosen to be $1{,}000$
iterations and the subsequent $1{,}000$ iterations were used to compute
the posterior mean. Increasing the burn-in period to more than $1{,}000$
iterations results in insignificant changes to our estimates, thus
indicating that the burn-in period chosen is sufficient. The posterior
mean estimates for both the priors, together with the MLE estimates,
are provided in Table~\ref{icfestmbys}. The running time for the Gibbs
sampling procedure for each prior is approximately $26$ seconds on a
Pentium M $1.6$ GHz processor. We find that the Bayesian approach using
our priors and the corresponding block Gibbs sampler gives stable
estimates and thus yields a useful alternative methodology for
inference in covariance graph models. We also note that the two
different vertex orderings yield very similar results.

\begin{table}
\tabcolsep=0pt
\caption{ICF: Maximum likelihood estimate from iterative conditional
fitting. \textup{BY1}: Bayesian posterior mean estimate for prior $\pi_{U_1,
\bolds{\alpha}^1}$. \textup{BY2}: Bayesian posterior mean\vspace*{-1pt} estimate for prior
$\pi_{U_2,\bolds{\alpha}^2}$. $\mathrm{\widetilde{BY1}}$: Bayesian posterior
mean estimate for prior $\pi_{U_1, \bolds{\alpha}^1}$ with\vspace*{-1pt} a different
ordering. $\mathrm{\widetilde{BY2}}$: Bayesian posterior mean estimate for prior
$\pi_{U_2,\bolds{\alpha}^2}$ with a different ordering\vspace*{-2pt}}\label{icfestmbys}
\begin{tabular*}{\textwidth}{@{\extracolsep{\fill}}lccccccccc@{}}
\hline
 & \textbf{GAL11} & \textbf{GAL4} & \textbf{GAL80} & \textbf{GAL3} & \textbf{GAL7} & \textbf{GAL10}
 & \textbf{GAL1} & \textbf{GAL2} & \textbf{Method} \\
\hline
GAL11 & 0.152 & 0.030 & 0\phantom{000,} & $-$0.052\phantom{$-$} & 0\phantom{000,} & 0\phantom{000,} & 0\phantom{000,} & $-$0.068\phantom{$-$} & ICF\\
 & 0.164 & 0.030 & 0\phantom{000,} & $-$0.050\phantom{$-$} & 0\phantom{000,} & 0\phantom{000,} & 0\phantom{000,} & $-$0.068\phantom{$-$} & BY1\\
 & 0.156 & 0.030 & 0\phantom{000,} & $-$0.052\phantom{$-$} & 0\phantom{000,} & 0\phantom{000,} & 0\phantom{000,} & $-$0.068\phantom{$-$} & BY2\\
 & 0.152 & 0.030 & 0\phantom{000,} & $-$0.051\phantom{$-$} & 0\phantom{000,} & 0\phantom{000,} & 0\phantom{000,} & $-$0.069\phantom{$-$} & $\widetilde{\mathrm{BY1}}$\\
 & 0.155 & 0.030 & 0\phantom{000,} & $-$0.052\phantom{$-$} & 0\phantom{000,} & 0\phantom{000,} & 0\phantom{000,} & $-$0.070\phantom{$-$} & $\widetilde{\mathrm{BY2}}$\\
 [3 pt]
GAL4 & & 0.128 & 0.040 & 0.042 & 0\phantom{000,} & 0\phantom{000,} & 0\phantom{000,} & 0.030 & ICF\\
 & & 0.142 & 0.040 & 0.041 & 0\phantom{000,} & 0\phantom{000,} & 0\phantom{000,} & 0.027 & BY1\\
 & & 0.133 & 0.041 & 0.042 & 0\phantom{000,} & 0\phantom{000,} & 0\phantom{000,} & 0.028 & BY2\\
 & & 0.128 & 0.039 & 0.040 & 0\phantom{000,} & 0\phantom{000,} & 0\phantom{000,} & 0.026 & $\widetilde{\mathrm{BY1}}$\\
 & & 0.132 & 0.040 & 0.042 & 0\phantom{000,} & 0\phantom{000,} & 0\phantom{000,} & \phantom{0}0.0278 & $\widetilde{\mathrm{BY2}}$\\
 [3 pt]
GAL80 & & & 0.223 & 0.082 & 0.197 & 0.198 & 0.239 & 0.227 & ICF\\
 & & & 0.237 & 0.072 & 0.193 & 0.194 & 0.235 & 0.216 & BY1\\
 & & & 0.232 & 0.076 & 0.199 & 0.2\phantom{00} & 0.243 & 0.223 & BY2\\
 & & & 0.224 & 0.076 & 0.197 & 0.197 & 0.240 & 0.218 & $\widetilde{\mathrm{BY1}}$\\
 & & & 0.232 & 0.076 & 0.202 & 0.203 & 0.245 & 0.227 & $\widetilde{\mathrm{BY2}}$\\ [3pt]
GAL3 & & & & 0.612 & 0.723 & 0.549 & 0.515 & 0.582 & ICF\\
 & & & & 0.626 & 0.713 & 0.544 & 0.509 & 0.575 & BY1\\
 & & & & 0.643 & 0.747 & 0.568 & 0.532 & 0.599 & BY2\\
 & & & & 0.628 & 0.719 & 0.549 & 0.517 & 0.582 & $\widetilde{\mathrm{BY1}}$\\
 & & & & 0.667 & 0.749 & 0.574 & 0.531 & 0.605 & $\widetilde{\mathrm{BY2}}$\\ [3 pt]
GAL7 & & & & & 3.422 & 2.593 & 2.768 & 2.540 & ICF\\
 & & & & & 3.462 & 2.584 & 2.756 & 2.533 & BY1\\
 & & & & & 3.588 & 2.682 & 2.866 & 2.636 & BY2\\
 & & & & & 3.541 & 2.588 & 2.761 & 2.532 & $\widetilde{\mathrm{BY1}}$\\
 & & & & & 3.708 & 2.681 & 2.865 & 2.627 & $\widetilde{\mathrm{BY2}}$\\ [3 pt]
GAL10 & & & & & & 2.372 & 2.409 & 2.267 & ICF\\
 & & & & & & 2.373 & 2.400 & 2.266 & BY1\\
 & & & & & & 2.453 & 2.497 & 2.358 & BY2\\
 & & & & & & 2.389 & 2.407 & 2.277 & $\widetilde{\mathrm{BY1}}$\\
 & & & & & & 2.473 & 2.489 & 2.356 & $\widetilde{\mathrm{BY2}}$\\ [3 pt]
GAL1 & & & & & & & 2.890 & 2.502 & ICF\\
 & & & & & & & 2.961 & 2.501 & BY1\\
 & & & & & & & 3.086 & 2.604 & BY2\\
 & & & & & & & 2.969 & 2.496 & $\widetilde{\mathrm{BY1}}$\\
 & & & & & & & 3.087 & 2.582 & $\widetilde{\mathrm{BY2}}$\\ [3pt]
GAL2 & & & & & & & & 2.870 & ICF\\
 & & & & & & & & 3.003 & BY1\\
 & & & & & & & & 3.153 & BY2\\
 & & & & & & & & 2.892 & $\widetilde{\mathrm{BY1}}$\\
 & & & & & & & & 3.005 & $\widetilde{\mathrm{BY2}}$\\
\hline
\end{tabular*}
\end{table}

\subsection{Simulation example} \label{blkgssimel}

\begin{figure}

\includegraphics{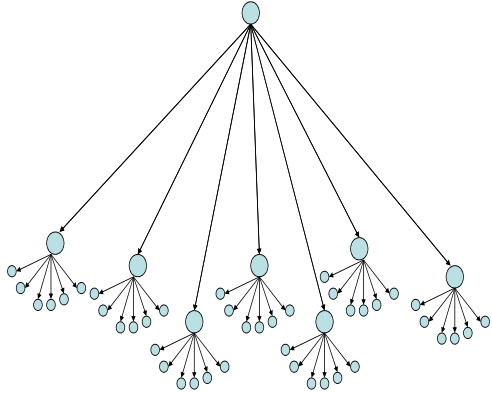}

\caption{Hasse diagram for a homogeneous graph with
$50$ vertices.} \label{exmlhgsfigu}
\end{figure}

A proof of convergence of the block Gibbs sampling algorithm proposed
in Section \ref{blkgssimul} was provided in Section \ref{gibbscnvgc}.
The speed at which convergence occurs is also a very important concern
for implementation of the algorithm. The number of steps that are
required before one can generate a reasonable approximate sample from
the posterior distribution is reflective of the rate of convergence.
Understanding this is important for the accuracy of Bayes estimates
such as   the posterior mean. We proceed to investigate the performance
of the block Gibbs sampling algorithm in a situation where the
posterior mean is known exactly and hence allows a direct comparison.
Consider a homogeneous graph $G$ with $50$ vertices, with the
corresponding Hasse tree given by Figure~\ref{exmlhgsfigu}. Let $\Sigma
\in P_G$, where the vertices have been ordered according to the Hasse
perfect vertex elimination scheme of Section \ref{hmgnsgphpy}, the
diagonal entries are $50$ and all other nonzero entries are $1$. We
simulate $100$ observation vectors ${\mathbf Y}_1, {\mathbf Y}_2,
{\mathbf Y}_3, \ldots, {\mathbf Y}_{99}, {\mathbf Y}_{100}$ from
$\mathcal{N}_{50} (0, \Sigma)$. For illustration purposes, we choose a
diffuse prior $\pi^{\Sigma}_{U,\bolds{\alpha}}$ with $U = 0$ and
$\alpha_i = 2 |\mathcal{N}^< (i)| + 5,  i = 1,2, \ldots, 50$.

Since the graph $G$ is homogeneous, we can compute the
posterior mean $\Sigma_{\mathrm{mean}} := {\mathbf E}_{U,\bolds{\alpha}}[
\Sigma \mid {\mathbf Y}_1, {\mathbf Y}_2, \ldots, {\mathbf Y}_{100}]$
explicitly. We can therefore assess the ability of the block Gibbs
sampling algorithm to estimate the posterior mean by comparing it to
the true value of the mean. We run the block Gibbs sampling algorithm
to sample from the posterior distribution and subsequently check its
performance in estimating $\Sigma_{\mathrm{mean}}$. We use an initial burn-in
period of $B$ iterations and then average over the next $I$ iterations
to get the estimate $\widehat{\Sigma}$. The times needed for
computation (using the R software) and the\vspace*{1pt} relative errors $\frac{\|
\widehat{\Sigma} - \Sigma_{\mathrm{mean}} \|_2}{\| \Sigma_{\mathrm{mean}} \|_2}$
corresponding to various choices of $B$ and $I$ are provided in Table
\ref{diagnostic}. The diagnostics in Table~\ref{diagnostic} indicate
that the block Gibbs sampling algorithm performs exceptionally well,
yielding estimates that approach the true mean in only a few thousand
steps. The time taken for running the algorithm is also provided in
Table \ref{diagnostic}.

The diagnostics in Table \ref{diagnostic} indicate that the block Gibbs
sampling algorithm performs exceptionally well, yielding estimates that
approach the true mean in only a few thousand steps. The time taken for
running the algorithm is also provided in Table \ref{diagnostic}.

\begin{table}[b]\tablewidth=250pt
\caption{Performance assessment of the Gibbs sampler in simulation
example} \label{diagnostic}
\begin{tabular*}{250pt}{@{\extracolsep{\fill}}lccc@{}}
\hline
\textbf{Burn-in (B)} & \textbf{Average (I)} & \textbf{Time (seconds)} & \textbf{Relative error}\\
\hline
1000 & 1000 & 139.77 & 0.01748220\\
2000 & 1000 & 209.72 & 0.01240595\\
3000 & 1000 & 279.52 & 0.01300910\\
4000 & 1000 & 349.44 & 0.01142864\\
4000 & 3000 & 489.19 & 0.01246141\\
4000 & 5000 & 631.21 & 0.01081264\\
4000 & 7000 & 769.70 & \phantom{0}0.009244206\\
\hline
\end{tabular*}
\end{table}

\section{Closing remarks}

In this paper, we have proposed a theoretical framework for Bayesian
inference in covariance graph models. The main challenge was the
unexplored terrain of working with curved exponential families in the
continuous setting. A rich class of conjugate priors has been developed
in this paper for covariance graph models where the underlying graph is
decomposable.

We have been able to exploit the structure of the conjugate priors to
develop a block Gibbs sampler to effectively sample from the posterior
distribution. A rigorous proof of convergence is also given. Comparison
with other classes of priors is also undertaken. We are able to compute
the normalizing constant for homogeneous graphs, thereby making
Bayesian model selection possible in a tractable way for this class of
models. The Bayesian approach yields additional dividends, in the sense
that we can now carry out inference in covariance graph models, even
when the sample size $n$ is less than the dimension $p$ of the data,
something which is otherwise not generally possible in the maximum
likelihood framework. Furthermore, we thoroughly explore the
theoretical properties of our class of conjugate priors. In particular,
in the homogeneous case, hyper-Markov properties and closed form
expressions for the expected value of the covariance matrix are
established.  Furthermore, the usefulness of the methodology that is
developed is illustrated through examples. A couple of open problems
are worth mentioning:
\begin{itemize}
\item What are the necessary conditions for the existence of the
normalizing constant for decomposable graphs?
\item Does the
hyper-Markov property for the class of priors developed in this paper
hold for decomposable graphs?
\end{itemize}
We conclude by noting that the use of the class of Wishart
distributions introduced in this paper for Bayesian inference, along
with a detailed study of Bayes estimators in this context, is clearly
an important topic and is the focus of current research.

\begin{appendix}\label{app}
\section*{Appendix}

\begin{pf*}{Proof of Proposition \ref{inverselrt}}
From the definition of $N$
and $L$, it is easy to verify that
\begin{eqnarray*}
& & (LN)_{ii} = 1\qquad  \forall 1 \leq i \leq m,\\
& & (LN)_{ij} = 0\qquad  \forall 1 \leq i < j \leq m.
\end{eqnarray*}
Now, let $i > j$. It follows, by the definition of $N,$ that
\begin{eqnarray*}
(LN)_{ij}&=& \sum_{k=j}^i L_{ik} N_{kj}\\
&=& N_{ij} + \sum_{k=j+1}^{i-1} L_{ik} \sum_{\bolds{\tau} \in A,
\tau_1 = k, \tau_{\operatorname{dim}(\bolds{\tau})} = j} (-1)^{\operatorname{dim}(\bolds{\tau})-1} L_{\bolds{\tau }}+ L_{ij}\\
&=& N_{ij} - \sum_{k=j+1}^{i-1}  \sum_{\bolds{\tau} \in A, \tau_1 = k,
\tau_{\operatorname{dim}(\bolds{\tau})} = j} (-1)^{\operatorname{dim}(\bolds{\tau})} L_{ik}L_{\bolds{\tau}} + L_{ij}.
\end{eqnarray*}
Note that any $\bolds{\tau}' \in \mathcal{A}$ with $ \bolds{\tau}'_1 =
i$, $\bolds{\tau}'_{\operatorname{dim}(\bolds{\tau}')} = j$, $\operatorname{dim}(\bolds{\tau}')
> 2$ can be uniquely expressed as $\bolds{\tau}' = (i, \bolds{\tau})$,
where $j+1 \leq \bolds{\tau}_1 \leq i-1,  \tau_{\operatorname{dim}(\bolds{\tau})} =
j$. Recall that, by definition, $L_{\bolds{\tau}'} = L_{i \tau_1}
L_{\bolds{\tau}}$. Also, if $ \bolds{\tau}' \in \mathcal{A}$ with
$\bolds{\tau}'_1 = i$, $\bolds{\tau}'_{\operatorname{dim}(\bolds{\tau}')} = j$,
$\operatorname{dim}(\bolds{\tau}') = 2$, then $\bolds{\tau}' = (i,j)$ and
$L_{\bolds{\tau}'} = L_{ij}$. Hence,
\begin{eqnarray*}
(LN)_{ij} &=& N_{ij} - \sum_{\bolds{\tau}' \in A, \tau'_1 = i,
\tau'_{\operatorname{dim}( \bolds{\tau}')} = j}
(-1)^{\operatorname{dim}(\bolds{\tau}')-1} L_{\bolds{\tau}'}\\
&=& N_{ij} - N_{ij}\\
&=& 0.
\end{eqnarray*}
Hence, $LN = I$ and thus $L^{-1} = N$.
\end{pf*}

\begin{pf*}{Proof of Theorem \ref{nrmcnstdml}} Let us simplify the integral
by integrating out the terms $D_{ii}$, $1 \leq i \leq m$:
%
\begin{eqnarray}
& & \int e^{- ( \operatorname{tr} ( (LDL^T)^{-1} U ) + \sum_{i=1}^m \alpha_i
\log D_{ii})/2}\,dL \,dD\nonumber\\
&& \qquad=\int e^{- ( \operatorname{tr} ( D^{-1} ( L^{-1} U (L^T)^{-1} ) ) +
\sum_{i=1}^m \alpha_i \log D_{ii} )/2}\, dL\, dD\nonumber\\
&& \qquad=\int \prod_{i=1}^m e^{- ( L^{-1} U (L^T)^{-1})_{ii}/(2D_{ii})}
D_{ii}^{-\alpha_i/2}\,dD\, dL\nonumber\\
&&\qquad= \int \prod_{i=1}^m \frac{\Gamma ( \alpha_i/2 - 1 )
2^{\alpha_i/2 - 1}}{( ( L^{-1} U (L^T)^{-1} )_{ii})^{\alpha_i/2 - 1}}\,dL\nonumber\\
\eqntext{(\mbox{assuming } \alpha_i > 2\ \forall
i =1, 2, \ldots, m)}\\
&&\qquad= \int \prod_{i=1}^m \frac{\Gamma ( \alpha_i/2 - 1 )
2^{\alpha_i/2 - 1}}{( (L^{-1})_{i \cdot} U ((L^{-1})_{i
\cdot})^T )^{\alpha_i/2 - 1}}\,dL. \qquad(**)\nonumber
\end{eqnarray}
In order to simplify this integral, we perform a change of measure by
transforming the nonzero elements of $L$ to the corresponding elements
of $L^{-1}$. For convenience and brevity, the notation $L^{-1}_{ij}$ is
used in place of $(L^{-1})_{ij}$. Now, note the following facts.
\begin{enumerate}
\item Let $L \in \mathcal{L}_G$. From Proposition \ref{inverselrt}, for
$(i,j) \in E$, $i > j$,
%
\begin{equation} \label{invrsbjctn}
L^{-1}_{ij} = - L_{ij} + f \bigl( (L_{uv})_{(u,v) \in E,  j \leq u < i,  j
\leq v < u\  \mathrm{or}\  u = i, j < v < i} \bigr),
\end{equation}
that is, $L^{-1}_{ij} + L_{ij}$ is a function ($f$) of $L_{uv}, (u,v)
\in E,  j \leq u < i,  j \leq v < u \mbox{ or } u = i$, $j < v < i$,
such that $f$ is zero when all its arguments are zero.
We use the above to show that $L$ is a function of
$\{L^{-1}_{uv}\}_{u>v,(u,v) \in E}$. Let $i^* = \min \{i\dvtx L_{ij} \neq 0
\mbox{ for some } j < i \}$. Let $j^* = \max \{j\dvtx L_{i^* j} \neq 0\}$.
By (\ref{invrsbjctn}) and the definition of $i^*$ and $j^*$, we have
$L_{i^* j^*} = - L^{-1}_{i^* j^*}$. We proceed by induction. Let $i >
j$, $(i,j) \in E$ and suppose that the hypothesis is true for all
$(u,v) \in E,  1 \leq u < i$, $1 \leq v < u$ or $u = i$, $j < v < i$.
Then,
\[
L_{ij} = - L^{-1}_{ij} + f \bigl( (L_{uv})_{(u,v) \in E,  j \leq u < i,  j
\leq v < u\  \mathrm{or}\  u = i,  j < v < i} \bigr)
\]
and, by the induction hypothesis, the right-hand side of the above
equation is a function of $\{L^{-1}_{uv}\}_{u>v, (u,v) \in E}$. Hence,
the matrix $L$ is a function of $\{L^{-1}_{uv}\}_{u>v,(u,v) \in E}$.

It follows that the transformation
\[
\{L_{ij}\}_{(i,j) \in E, i > j} \rightarrow \{L^{-1}_{ij}\}_{(i,j) \in
E, i > j}
\]
is a bijection and the absolute value of the Jacobian of this
transformation is $1$ since  it is the determinant of a
lower-triangular matrix with diagonal entries $1$.
\item If ${\mathbf x} = \left({{\mathbf{x}_1}\atop{\mathbf{x}_2}}\right)$
and $U =\left({{U_{11} \enskip U_{12}}\atop{U_{21} \enskip U_{22}}}\right)$
is a positive definite matrix, then
%
\begin{equation} \label{pdmatrixss}
\quad{\mathbf x}^T U {\mathbf x} = z^T z + {\mathbf x}_2^T ( U_{22} - U_{21}
U_{11}^{-1} U_{12} ) {\mathbf x}_2 \geq {\mathbf x}_2^T ( U_{22} -
U_{21} U_{11}^{-1} U_{12} ) {\mathbf x}_2,
\end{equation}
where $z = U_{11}^{1/2} {\mathbf x}_1 + U_{11}^{-1/2}
U_{12} {\mathbf x}_2$.
\end{enumerate}
Hence, after transforming the nonzero entries of $L$ to the
corresponding entries of $L^{-1}$ and using (\ref{pdmatrixss}) to
eliminate the dependent entries of $L^{-1}$ from the integrand, we get
\begin{eqnarray*}
& & \int e^{- ( \operatorname{tr}( (LDL^T)^{-1} U ) + \sum_{i=1}^m \alpha_i
\log D_{ii} )/2}\,dL\,dD\\
&&\qquad= \int \prod_{i=1}^m \frac{\Gamma (\alpha_i/2 - 1 )
2^{\alpha_i/2 - 1}}{( (L^{-1})_{i \cdot} U ((L^{-1})_{i\cdot})^T)^{\alpha_i/2 - 1}}\,dL\\
&&\qquad\leq K \prod_{i=2}^m \int_{\mathbb{R}^{|\mathcal{N}^\prec (i)|}}
\frac{1}{\left(\pmatrix{
{\mathbf a}_i^T & 1} U_i^* {{\mathbf a}_i\choose
 1 }\right)^{\alpha_i/2 - 1}}\,d {\mathbf a}_i.
\end{eqnarray*}
Here, $K$ is a constant, $U_i^*$ is an appropriate positive definite
matrix and ${\mathbf a}_i$ represents the independent entries in the
$i$th row of $L^{-1}$. By a suitable linear transformation ${\mathbf
b}_i$ of each of the ${\mathbf a}_i$, $i = 2,3, \ldots, m$, we get
\begin{eqnarray*}
& & \int e^{- ( \operatorname{tr} ( (LDL^T)^{-1} U ) + \sum_{i=1}^m
\alpha_i \log D_{ii} )/2}\, dL \,dD\\
&&\qquad\leq K^* \prod_{i=2}^m \int_{\mathbb{R}^{|\mathcal{N}^\prec (i)|}}
\frac{1}{( {\mathbf b}^T {\mathbf b} + u_i^{**} )^{\alpha_i/2 -1}}\,d {\mathbf b}_i.
\end{eqnarray*}
Here, $K^*$ and $u_i^{**}$, $i = 2,3, \ldots, m,$ are constants. Using
the standard fact that
\[
\int_{\mathbb{R}^k} \frac{1}{( {\mathbf x}^T {\mathbf x} + 1
)^\gamma}\,
d {\mathbf x} < \infty \qquad\mbox{if } \gamma > \frac{k}{2},
\]
we conclude that
\[
\int e^{- ( \operatorname{tr}( (LDL^T)^{-1} U ) + \sum_{i=1}^m \alpha_i \log
D_{ii} )/2}\,dL\, dD < \infty
\]
if $\alpha_i > |\mathcal{N}^\prec (i)| + 2$ for all $i = 1,2, \ldots,
m$.
\end{pf*}
\end{appendix}

\section*{Acknowledgments}
The authors gratefully acknowledge the faculty at the Stanford
Statistics Department for their feedback and tremendous enthusiasm for
this work. We thank Persi Diaconis in particular for taking the time to
read the paper and for his encouraging remarks. The second-named author
gratefully acknowledges Professor Phil Dawid, Professor Steffen
Lauritzen and the London Mathematical Society (LMS) for the kind
invitation to attend the LMS-sponsored meetings in Durham in June/July
2008. The author also gratefully acknowledges Professor Nanny Wermuth
for an invitation to a graphical models workshop in Sweden in July
2008. In addition, the Isaac Newton Institute for Mathematical Sciences
at Cambridge University, U.K., is gratefully acknowledged for generous
support and for use of their facilities, where the second named author
was in residence for part of the 2008 program on ``Statistical Theory
and Methods for Complex, High-Dimensional Data.'' The organizers of
this program and a travel award from the Cambridge Philosophical
Society are also acknowledged. Bala Rajaratnam gratefully also
acknowledges support from an NSERC fellowship. The authors would like
to thank organizers of the O-Bayes meeting in Philadelphia in June,
2009, and Professor Lawrence Brown in particular for very positive
feedback and encouraging words. The authors also acknowledge two
anonymous referees for useful comments.

\printaddresses

\end{document}